\documentclass[moor]{informs1}              

\usepackage{natbib}
 \NatBibNumeric
 \def\bibfont{\small}%
 \bibpunct[, ]{[}{]}{,}{n}{}{,}%

\usepackage[colorlinks=true,breaklinks=true,bookmarks=true,urlcolor=blue,
     citecolor=blue,linkcolor=blue,bookmarksopen=false,draft=false]{hyperref}

\def\EMAIL#1{\href{mailto:#1}{#1}}
\def\URL#1{\href{#1}{#1}}         

\TheoremsNumberedBySection  

\EquationsNumberedBySection 

\newcommand{\opt}{\mathrm{OPT}}
\newcommand{\mb}[1]{\ensuremath{\boldsymbol{#1}}}
\usepackage{amsmath,  amssymb, epsfig, graphicx, float, url,tabularx,multirow}
\usepackage{algorithm, algorithmicx, algpseudocode,bbm,subcaption,booktabs}
\usepackage{enumitem}
\usepackage{footnote}
\makesavenoteenv{tabular}
\makesavenoteenv{table}

\newcommand*{\red}{\textcolor{black}}

\begin{document}



\RUNAUTHOR{El Housni and Goyal}
 \RUNTITLE{On the Optimality of Affine Policies for Budgeted Uncertainty Sets}

\TITLE{On the Optimality of Affine Policies for Budgeted Uncertainty Sets}
\ARTICLEAUTHORS{%
\AUTHOR{Omar El Housni}
\AFF{Dept of Industrial Engineering and Operations Research, Columbia University, New York NY 10027, \EMAIL{oe2148@columbia.edu}, \URL{http://www.columbia.edu/~oe2148/}}
\AUTHOR{Vineet Goyal}
\AFF{Dept of Industrial Engineering and Operations Research, Columbia University, New York NY 10027, \EMAIL{vg2277@columbia.edu}, \URL{http://www.columbia.edu/~vg2277/}}
} 

\ABSTRACT{%
In this paper, we study the performance of affine policies for two-stage adjustable robust optimization problem with \red{fixed recourse} and uncertain right hand side belonging to a budgeted uncertainty set. This is an important class of uncertainty sets widely used in practice where we can specify a budget on the adversarial deviations of the uncertain parameters from the nominal values to adjust the level of conservatism. The two-stage adjustable robust optimization problem is hard to approximate within a factor better than $\Omega( \frac{\log n}{\log \log n})$   even for budget of uncertainty sets 
where $n$ is the number of decision variables. Affine policies, where the second-stage decisions are  constrained to be an affine function of the uncertain parameters, provide a tractable approximation for the problem \red{and have been observed to exhibit good empirical performance}. \red{We show that affine policies  give an $O( \frac{\log n}{\log \log n})$-approximation for the two-stage adjustable robust problem with fixed non-negative recourse for budgeted uncertainty sets.} This matches the hardness of approximation and therefore, surprisingly affine policies provide an optimal approximation for the problem (up to a constant factor). We also show strong theoretical performance bounds for affine policy for significantly more general class of intersection of budgeted sets including disjoint constrained budgeted sets, permutation invariant sets and general intersection of budgeted sets. Our analysis relies on showing the existence of a near-optimal feasible affine policy that satisfies certain nice structural properties. Based on these structural properties, we also present an alternate algorithm to compute near-optimal affine solution that is significantly faster than computing the optimal affine policy by solving a large linear program.  
}%


\MSCCLASS{Primary: 90C15 , 90C47}

\maketitle

\section{Introduction.}
             
In this paper, we consider the following two-stage adjustable robust linear optimization problem with uncertain right hand side:
\begin{equation}\label{eq:ar}
\begin{aligned}
z_{\sf AR}({\cal U}) = \min_{\mb x} \; & \mb{c}^T \mb{x} + \max_{\mb{h}\in {\cal U}} \min_{\mb{y}(\mb{h})} \mb{d}^T \mb{y}(\mb{h}) \\
& \mb{A}\mb{x} + \mb{B}\mb{y}(\mb{h}) \; \geq \; \mb{h},   \; \; \; \forall \mb h \in {\cal U} \\
& \red{\mb{x}  \in  {\cal X} }\\
& \mb{y}(\mb{h}) \in  {\mathbb R}^{n}_+, \;  \; \; \forall \mb h \in {\cal U},
\end{aligned}
\end{equation}
where $\mb{A} \in {\mathbb R}^{m {\times} n},\mb{B} \in {\mathbb R}^{m {\times} n}_+, \mb{c}\in {\mathbb R}^{n}_+$ and $ \mb{d}\in {\mathbb R}^{n}_+$ . The right-hand-side $\mb{h}$ is uncertain and belongs to a compact convex uncertainty set ${\cal U}\subseteq{\mathbb R}^m_+$. \red{ The recourse matrix is non-negative and fixed, i.e., $\mb B$ belongs to the non-negative orthant and does not depend on the uncertain parameter $\mb h$.}
 \red{The goal in this problem is to select the first-stage decision $\mb{x} \in {\cal X}$, where $ {\cal X}$ is a polyhedral set }and the second-stage recourse decision,  $\mb y(\mb{h})$, as a function of the uncertain right hand side realization, $\mb{h}$ such that the worst-case cost over all realizations of $\mb h \in {\cal U}$ is minimized. We assume that  ${\cal U} \subseteq [ 0, 1]^m$ and  $\forall i \in [m], \mb e_i \in {\cal U}$. This assumption is without loss of generality since we can scale the constraint matrices $\mb A$ and $\mb B$ to satisfy the assumption without changing the optimal.

This model has been widely considered in the literature (see for example 
Bertsimas and de Ruiter \cite{bertsimas2016duality}, Bertsimas and Goyal \cite{BG10}, Dhamdhere et al. \cite{dhamdhere2005pay}, El Housni and Goyal \cite{NIPS2017_7061}, Gupta et al. \cite{gupta2014thresholded}, Xu and Burer \cite{Xu2018}, Zhen et al. \cite{zhen2016adjustable}.)
It captures many important applications including set cover, capacity planning and network design problems under uncertain demand. Here, the right hand side $\mb h$ models the uncertain demand and the covering constraints capture the requirement of satisfying the uncertain demand. However, the adjustable robust optimization problem \eqref{eq:ar} is intractable in general. In fact, Feige et al.~\cite{FJMM07} show that the two-stage adjustable  problem \eqref{eq:ar} can not be approximated within a ratio better than $\Omega(\frac{\log n}{\log \log n})$  under a reasonable complexity assumption, namely, 3SAT can not be solved in time $2^{O(\sqrt{n})}$ on instances of size $n$. 

In view of the intractability, several approximation policies (or \textit{decision rules}) have been considered in the literature for \eqref{eq:ar} including static, piecewise static, affine and piecewise affine policies. In a static policy, we compute a single optimal solution $(\mb{x},\mb{y})$ that is feasible for all realizations of the uncertain right hand side. Bertsimas et al.~\cite{BGS10}  relate the performance of static solution to the symmetry of the uncertainty set and show that it provides a good approximation to the adjustable problem if the uncertainty set verifies some symmetry properties. However, static policy is too conservative in general and the performance of static solutions can be arbitrarily large for a general convex uncertainty set. 

Ben-Tal et al. \cite{Ben-Tal04} introduce affine policy approximation for \eqref{eq:ar}, where they restrict the second-stage decision, $\mb{y}(\mb{h})$ to being an affine function of the uncertain right-hand-side $\mb{h}$, i.e.,  $\mb{y}(\mb{h})=\mb{P}\mb{h}+\mb{q}$ for some decision variables $\mb{P}\in{\mathbb R}^{n\times m}$ and $\mb{q}\in{\mathbb R}^n$. Affine policy can be computed efficiently for a large class of uncertainty sets and therefore, provide a tractable approximation for the two-stage problem. Furthermore, the empirical performance of affine policies has been observed to be near-optimal for a large class of instances even though theoretically, optimality of affine policies is known in only a few settings. Bertsimas et al. \cite{BIP09} and Iancu et al. \cite{ISS13} show that affine policy is optimal for multi-stage adjustable problems with a single uncertain parameter at each stage. Bertsimas and Goyal \cite{BG10} show that affine policy is optimal for the two-stage adjustable problem \eqref{eq:ar} if the uncertainty set ${\cal U}$ is a simplex. \red{However, in the particular case where we assume only non-negativity constraints on the first stage decision variable, i.e. $\mb x \geq \mb 0$, they show that
the worst-case performance of an optimal affine solution is $\Theta(\sqrt m)$ times the optimal cost of~\eqref{eq:ar}~\cite{BG10}. Note that the gap could be even larger for general polyhedral constraints that involves only $\mb x$ i.e., $\mb x \in {\cal X}$.} Therefore, there is a significant gap between the worst-case performance of affine policies and the observed empirical performance.

More general decision rules have been considered in the literature for two-stage problems; piecewise affine policies (Ben-Tal et al. \cite{ben2018tractable}), binary decision rules (Bertsimas and Georghiou \cite{bertsimas2015design}), adjustable solutions via iterative splitting of uncertainty sets (Postek and Den Hertog \cite{postek2016multistage}), k-adaptibility (Hanasusanto et al. \cite{hanasusanto2015k}), extended affine decision rules (Chen et al. \cite{chen2008linear}), etc. While these decision rules can improve in some instances over affine policies, they become computationally very challenging especially for large size instances. 
For a more extensive review of the  literature, we refer the reader to Bertsimas et al. \cite{BBC08} and Ben-Tal et al. \cite{BNE10}.

In this paper, we study the performance of affine policies for the two-stage adjustable problem \eqref{eq:ar}. One of our important goals in this work is to address the stark contrast between the observed near-optimal empirical performance and the worst-case approximation bound of $\Theta(\sqrt m)$~\cite{BG10}. Towards this, we consider the class of budget of uncertainty sets and intersection of budget of uncertainty sets that was introduced in Bertsimas and Sim~\cite{BS04}. This is widely used class of uncertainty sets in practice where the decision maker can specify a budget on the sum of adversarial deviations of the uncertain parameter from the nominal values. In particular, a budget of uncertainty set can be formulated as follows:
\begin{equation}\label{eq:budget}
{\cal U}= \left\{ \mb h \in {[ 0, 1]}^m \; \bigg\vert \; \sum _{i=1}^m  w_i h_i  \leq  1 \right\},
\end{equation}
where $ w_i \in [0,1]$ for all $i \in [m]$. It is known that the two-stage adjustable problem \eqref{eq:ar}  is hard to approximate under this class of uncertainty set. In particular, Feige et al. \cite{FJMM07} show that the two-stage adjustable problem \eqref{eq:ar} where $\cal U$ is the budget of uncertainty set \eqref{eq:budget} is hard to approximate within a factor $\Omega( \frac{\log n}{\log \log n})$ even when all $w_i$ are equal, $\mb A,\mb B$ are 0-1 matrices \red{ and ${\cal X}= \mathbb{R}_+^n.$} The analysis in Bertsimas and Goyal \cite{BG10} shows that affine policy gives  $O( \sqrt{m})$-approximation to the adjustable problem \eqref{eq:ar} under this class of uncertainty sets \red{ and ${\cal X}= \mathbb{R}_+^n.$}  Bertsimas and Bidkhori \cite{bertsimas2014performance} provide a geometric bound  $O ( \frac{k^2+mk}{k^2+m} )$ in the special case  of \eqref{eq:budget}  where all $w_i=1/k$ \red{ and ${\cal X}= \mathbb{R}_+^n.$ This bound is} also $O(\sqrt m)$ in the worst-case for $k=\sqrt m$.

The above two-stage adjustable robust problem~\eqref{eq:ar} has also been considered in the context of combinatorial optimization problems such as  \red{network design under demand uncertainty} where the constraint matrices, $\mb A, \mb B \in \{0,1\}^{m \times n}$ and the first-stage and second-stage decisions are constrained to be binary (see for instance, Dhamdhere et al.~\cite{dhamdhere2005pay}, Feige et al.~\cite{FJMM07}, Gupta et al. \cite{gupta2014thresholded} and \cite{gupta2016robust}). Feige et al.~\cite{FJMM07} and Gupta et al.~\cite{gupta2014thresholded} give an $O(\log n)$-approximation for~\eqref{eq:ar} for the special case when $\mb A = \mb B \in \{0,1\}^{m \times n}$, first and second-stage costs are proportional, i.e., $\mb d = \lambda \cdot \mb c$ for some constant $\lambda \geq 1$ and a budget of uncertainty set with $w_i = 1/k$. Gupta et al.~\cite{gupta2016robust} also consider more general uncertainty set, namely, intersection of  {\em p-system} and {\em q-knapsack}\footnote{p-system and  q-knapsack are generalizations of matroid constraints. For more details, we refer the reader to \cite{gupta2016robust}} and give $O(pq\log n)$-approximation for the two-stage problem. However, we would like to note that the focus of this stream of work is to design approximation algorithms for combinatorial optimization problems where the decisions are constrained to be binary. Moreover, the algorithms are not restricted to and do not necessarily give decision rules or functional policy approximations for the two-stage problem. 

In contrast, the focus of our work is to analyze the performance of affine policies for the two-stage adjustable robust problem~\eqref{eq:ar} that are widely used in practice and exhibit strong empirical performance. 

 \subsection{Our contributions.} Our main contributions are the following.

\begin{enumerate}
\item[(a)] {\bf Optimal Approximation for Budget of Uncertainty Sets.} We show that affine policy gives an optimal approximation for the two-stage adjustable robust problem for budgeted uncertainty sets. In particular, affine policy gives an $O(\frac{\log n}{\log \log n})$-approximation to the two-stage adjustable problem~\eqref{eq:ar} where $\cal U$ is a budget of uncertainty set~\eqref{eq:budget}.  This performance bound matches the hardness of approximation~\cite{FJMM07}; thereby, showing that surprisingly affine policies give an optimal approximation (up to some constant factor) for~\eqref{eq:ar} for budget of uncertainty sets. In other words, there is no polynomial time algorithm with worst-case approximation guarantee better than affine policies by more than some constant factor. This bound significantly improves over the previous known bound of $O(\sqrt{m})$ \cite{BG10, bertsimas2014performance} for budget of uncertainty sets. \red{ Moreover, our result holds for general polyhedral constraints on the first stage variable $ \mb x \in \cal X$. In particular, we can model for example upper bounds on $\mb x$, this is in contrast  with the previous bounds in the literature that have been shown only in the special case of  $\cal X = \mathbb{R}_+^n$. }

Our analysis relies on constructing a feasible affine solution whose worst-case cost is within a factor $O( \frac{\log n}{\log \log n} )$ of the optimal cost. In particular, we partition the components of $\cal U$ into {\em inexpensive} and {\em expensive} components based on a threshold and construct an affine solution that covers \red{only the inexpensive components using a linear solution}. \red{The remaining components are covered using a static solution.} We show that for an appropriately chosen threshold that depends on the optimal cost, such an affine solution gives an $O(\frac{\log n}{\log \log n})$-approximation for the two-stage problem for the budget of uncertainty set. 

Therefore, in addition to establishing the performance bound that matches the hardness of approximation, our analysis shows there is a near-optimal affine solution whose structure is closely related to {\em threshold policies} that are widely used in many applications. This structural property might be of independent interest and also gives an alternate faster algorithm for computing near-optimal affine solutions for budget of uncertainty sets as we discuss later.

\vspace{2mm}
\item [(b)]{\bf Intersection of Budgeted Sets.} 
We also consider a more general family of uncertainty sets, namely the following intersection of budget of uncertainty sets:
\begin{equation} \label{eq:intersection:set}
{\cal U} = \left\{ \mb h \in [0,1]^m  \; \Big\vert \;  \sum_{i \in S_{\ell}} w_{ \ell i} h_i \leq 1  \;\;  \forall \ell \in  [L] \right\},
\end{equation}
where $ \mb w_l \in [0,1]^m,$ and $S_{\ell} \subseteq [m]$ for all $\ell \in [L]$. The set  \eqref{eq:intersection:set} is defined by the intersection of $L$ budget constraints. These are an important generalization of the budget of uncertainty set \eqref{eq:budget} that are widely used in practice. They capture for instance CLT sets \cite{bandi2012tractable}
 and {\em inclusion-constrained budgeted} sets \cite{gounaris2014adaptive}. 
 \begin{enumerate}[label=(\roman*)]
 \item  We first consider the case when the family of subsets $S_{\ell}$ for $\ell \in [L]$ are disjoint. We refer to this class of sets as {\em disjoint constrained budgeted sets}. These are essentially Cartesian product of $L$ budget of uncertainty sets. We show that affine policy is near-optimal and gives an $O \left( \log^2 n/ \log \log n  \right)$-approximation to \eqref{eq:ar} for this class of sets. We would like to note that the bound is independent of $L$.  Similar to the case of budget of uncertainty sets, our analysis is based on constructing a near-optimal affine solution by partitioning components of ${\cal U}$ into {\em inexpensive} and {\em expensive} components using appropriate thresholds for each of the $L$ budgeted sets in the Cartesian product. However, in this case, we are able to relate the performance of our affine solution to only a lower bound of $z_{\sf AR}( {\cal U})$. In particular, we use an {\em online algorithm} for the {\em fractional covering problem} to both construct thresholds (and therefore, a feasible affine solution) as well as the lower bound of the optimal \red{value}. 
 
\item For general intersection of $L$ budgets. \red{Under the assumption that $\cal X$ is a polyhedral cone (for example $\cal X = \mathbb{R}_+^n$),}  we show that affine policy gives $O \left(  \log L \log n/ \log \log n  \right)$ to \eqref{eq:ar} for the case where $\cal U$ is {\em permutation invariant}.  We say that $\cal U$ is permutation invariant if for any $\mb h \in {\cal U}$ and any $\tau $ permutation of $[m]$, then $\mb h^{\tau}\in{\cal U}$ where $ h^\tau_i= h_{\tau(i)}$. This class captures many important sets such as CLT sets. The performance of affine policy depends on $L$ in this case but degrades gracefully. For general intersection of budgeted sets \red{and $\cal X$  a polyhedral cone}, we show  a worst-case bound of $O \left(  L \log n/ \log \log n  \right)$ on the performance of affine policy for \eqref{eq:ar}. We summarize our results in Table \ref{tab:results}.
 \end{enumerate}
 
 \vspace{2mm}
\item [(c)] {\bf Faster Algorithm to compute Near-Optimal Affine Solutions.} 
Based on the structural properties of the near-optimal affine policies constructed for analysis of performance,  we present an alternate algorithm to compute an approximate affine policy for \eqref{eq:ar} for budget of uncertainty sets that is significantly faster than computing optimal affine policy by solving a large LP. In particular, our construction partitions the components into {\em inexpensive} and {\em expensive} based on a threshold depending on the optimal cost and 
shows the existence of a near-optimal affine solution \red{ that covers a fraction of {\em inexpensive} components using a linear solution and the remaining components using a static solution. }

From an algorithmic perspective, while we do not know the optimal cost and therefore, the threshold, we can still use this structure of a near-optimal affine solution to construct a good affine solution. \red{ In particular, our approximate affine solution can be computed efficiently by solving a single LP covering problem with $O(n+m)$ second stage variables as opposed to $O(nm)$ second stage variables in the optimal affine policy. }
Our algorithm scales very well and is significantly faster than computing affine policies. \red{ For instance, for $m=n=100$, it takes several minutes to compute the optimal affine policy whereas our algorithm computes an approximate affine policy within a few seconds. The comparison becomes even more stark when we increase the problem size. In particular, for $m=n=200$, the average time for optimal affine policy is more than an hour, whereas our algorithm computes an approximate affine policy in under $2$ minutes.   Moreover, our solution remains within $15\%$ of the optimal affine solution and the sub-optimality gap does not increase with dimension in our numerical experiments.} We would like to note that since our approximate affine is based on the construction of affine policy in our analysis, the worst-case approximation bound for the faster algorithm is also $O (  \frac{\log n}{\log \log n} )$.

\red{
\item [(d)] {\bf General case of recourse matrix $\mb B$.} 
We show that the  assumption on the non-negativity of the recourse matrix $\mb B$ is crucial for obtaining any non-trivial theoretical bounds on the performance of affine policies. In particular, we give a family of instances of the two-stage adjustable problem where the recourse matrix $\mb B$ is a network matrix with entries in $\{-1,0,1\}$ and show that the gap between optimal affine and adjustable policies can be unbounded even for the single budget of uncertainty set. The second-stage matrix being a network matrix captures important applications including lot sizing and facility location. }
\end{enumerate}

\begin{table}[ht]
\centering
\begin{tabular} {|l|l|l|}
  \hline
& Uncertainty sets   & Our Bounds \\
 \hline
1.& Budget of uncertainty set  \eqref{eq:budget} & \large{$O\left( \frac{\log n}{\log \log n} \right)^*$}    \\
   \hline
2.& Disjoint Intersection of Budgeted sets \eqref{eq:partition-matroid} &  \large{$O \left( \frac{\log^2 n}{\log \log n}   \right)$}    \\
\hline
3. & Permutation Invariant Intersection of Budgeted sets  \eqref{eq:intersection:set}& \large{$O \left(  \frac{\log L \cdot \log n}{\log \log n}  \right)^{**}$} \\
\hline
4. & General Intersection of Budgeted sets   \eqref{eq:intersection:set}& \large{$O \left(  \frac{L \cdot \log n}{\log \log n}   \right)^{**}  $}  \\
\hline
\end{tabular}%
\caption{Our performance bounds for affine policy under different uncertainty sets including budget of uncertainty set and intersection of budgeted sets. $^*$ denotes that this bound mathches the approximation hardness of \eqref{eq:ar}. \red{$^{**}$ denotes that these bounds hold under the assumption that $\cal X$ is a polyhedral cone.} }\label{tab:results}
 \end{table}

\section{Affine policies: Preliminaries.}
Affine policies (also known as {\em linear decision rules}) are widely used in the literature of robust optimization. They were introduced by Ben-Tal et al. \cite{Ben-Tal04} for the two-stage adjustable problem \eqref{eq:ar}. In an affine solution, we restrict the second stage decision $\mb y( \mb h )$ to be an affine function of the uncertain parameter $\mb h $, i.e., $$\mb y( \mb h )= \mb P \mb h + \mb q,$$ and we optimize over the variables $\mb P$ and $\mb q$. The affine problem is formulated as:
\begin{equation}\label{eq:aff}
\begin{aligned}
z_{\sf Aff}({\cal U}) = \min_{\mb x, \mb P, \mb q} \; & \mb{c}^T \mb{x} + \max_{\mb{h}\in {\cal U}}  \mb{d}^T\left( \mb{P}\mb{h}+\mb{q}\right) \\
& \mb{A}\mb{x} + \mb{B}\left( \mb{P}\mb{h}+\mb{q}\right) \; \geq \; \mb{h},   \; \; \; \forall \mb h \in {\cal U} \\
& \mb{P}\mb{h}+\mb{q} \; \geq \; \mb{0},   \; \; \; \forall \mb h \in {\cal U} \\
& \red{\mb{x}  \in  {\cal X} }.\\
\end{aligned}
\end{equation}
Affine policy has been widely used as an approximation to \eqref{eq:ar} due to its tractability. In fact, Ben-Tal et al. \cite{Ben-Tal04} show that affine problems have an equivalent standard LP formulation when the uncertainty set is described by a polyhedral set. The size of the LP is polynomial in the size of the input parameters (i.e.,  number of  variables and constraints in $\eqref{eq:ar}$ and number of constraints in $\cal U$). For completeness, we briefly discuss the tractability and compact LP formulation of affine policies. Consider the following polyhedral uncertainty set
\begin{equation} \label{def:U}
{\cal U}=\{\mb{h}\in{\mathbb R}^m_+\;|\;    \mb R \mb h \leq \mb r \},
\end{equation}
where $ \mb R \in \mathbb{R}^{L \times m}$ and $\mb r \in \mathbb{R}^L$. The affine problem \eqref{eq:aff} can be reformulated as the following epigraph formulation:
 \begin{equation*}
\begin{aligned}
z_{\sf Aff}({\cal U}) = \min \; & \mb{c}^T \mb{x} + z\\
& z \geq   \mb{d}^T\left( \mb{P}\mb{h}+\mb{q}\right),  \; \; \; \forall \mb h \in {\cal U} \\
& \mb{A}\mb{x} + \mb{B}\left( \mb{P}\mb{h}+\mb{q}\right) \; \geq \; \mb{h},   \; \; \; \forall \mb h \in {\cal U} \\
& \mb{P}\mb{h}+\mb{q} \; \geq \; \mb{0},   \; \; \; \forall \mb h \in {\cal U} \\
&  \red{\mb{x}  \in  {\cal X} }, \; \mb P \in \mathbb{R}^{n \times m}, \; \mb q \in \mathbb{R}^n, \; z \in \mathbb{R}.   \\
\end{aligned}
\end{equation*}
Note that this formulation can have infinitely many constraints but the separation problem is tractable. For example, the separation problem for the first set of constraints is: Given $z, \mb x, \mb P, \mb q$ decide if
\begin{equation} \label{eq:separation}
z - \mb d^T \mb q \geq \max_{\mb h \geq \mb 0} \; \{ \mb{d}^T\mb{P}\mb{h}\; | \;  \mb R  \mb h \leq \mb r\}.
\end{equation}
This can be done efficiently by solving the above maximization LP. Moreover, Ben-Tal et al. \cite{Ben-Tal04} show that we can formulate \eqref{eq:aff} as a compact LP using standard techniques from duality. For instance, consider the first set of constraints \eqref{eq:separation}, by taking the dual of the maximization problem, the constraint becomes $$z - \mb d^T \mb q \geq \min_{\mb v \geq \mb 0} \;     \{\mb{r}^T\mb v \; | \; \mb R^T \mb v \geq \mb P^T \mb d\}.$$
We can then drop the min and introduce $\mb v $ as a variable. Hence, we obtain the following linear constraints: 
$$ z- \mb d^T \mb q \geq \mb r^T \mb v, \qquad \mb R^T \mb v \geq \mb P^T \mb d, \qquad \mb{v}  \geq \mb 0.$$ We can apply the same techniques for the other constraints.
For completeness, we restate the compact LP formulation of Ben-Tal et al. \cite{Ben-Tal04}  adapted to our problem in the following lemma. The proof of the lemma is deferred to Appendix \ref{apx-LP}.

\begin{lemma} \label{lem:LP-form}
The affine problem \eqref{eq:aff} can be formulated as the following LP
\begin{equation}\label{eq:aff:LP}
\begin{aligned}
z_{\sf Aff}({\cal  U}) =  \min \; &  \mb{c}^T \mb{x} + z \\
& z- \mb d^T \mb q \geq \mb r^T \mb v \\
& \mb R^T \mb v \geq \mb P^T \mb d \\
& \mb A \mb x + \mb B \mb q \geq \mb V^T \mb r\\
& \mb R^T \mb V \geq \mb I_m - \mb B \mb P \\
& \mb q \geq \mb U^T \mb r \\
& \mb U^T \mb R + \mb P \geq \mb 0 \\
& \red{\mb{x}  \in  {\cal X} }, \; \mb{v}   \in  {\mathbb R}^{L}_+, \; \mb{U}   \in  {\mathbb R}^{L \times n}_+, \; \mb{V}   \in  {\mathbb R}^{L \times m}_+ \\
&  \mb P \in \mathbb{R}^{n \times m}, \; \mb q \in \mathbb{R}^n, \; z \in \mathbb{R}.   \\
\end{aligned}
\end{equation}
\end{lemma}

\section{Performance analysis for budget of uncertainty sets.}

In this section, we consider the class of budget of uncertainty sets \eqref{eq:budget} given by 
\begin{equation*}
{\cal U}= \left\{ \mb h \in {[ 0, 1]}^m \; \bigg\vert \; \sum _{i=1}^m w_i  h_i  \leq  1 \right\}.
\end{equation*}
As we mention earlier, this class  is widely used in the literature of robust optimization both in theory and practice. It provides the flexibility to adjust the level of conservatism in terms of probabilistic bounds on constraint violations. A special case of this class is when  $w_i$ are all equal to $\frac{1}{k}$ for  some  parameter $k \in \mathbb{N}$. In particular, in this case we have
\begin{equation} \label{eq:k-ones}
{\cal U}= \left\{ \mb h \in {[ 0, 1]}^m \; \bigg\vert \; \sum _{i=1}^m  h_i  \leq  k \right\}.
\end{equation}
The parameter $k$ is the {\em budget of uncertainty} that controls the conservatism of the uncertainty model. This special class \eqref{eq:k-ones} of budgeted sets is also known as the {\em cardinality constrained set} or  {\em $k$-ones polytope}. 

The two-stage adjustable problem \eqref{eq:ar} is known to be  $\Omega( \frac{\log n}{\log \log n})$-hard to approximate under the class of budget of uncertainty sets even in the special case of \eqref{eq:k-ones} (Feige et al. \cite{FJMM07}). We show that surprisingly the performance bound for affine policy matches this hardness of approximation.  In particular, we show that affine policy gives $O (\frac{\log n}{\log \log n}) $-approximation for~\eqref{eq:ar} under budget of uncertainty sets \eqref{eq:budget}.

\begin{theorem}\label{thm:budget}
Consider the two-stage adjustable problem \eqref{eq:ar} where $ {\cal U}$ is the budget of uncertainty set \eqref{eq:budget} . Then,
$$ z_{\sf Aff}({\cal U}) = O \left( \frac{\log n}{\log \log n} \right)\cdot z_{\sf AR}({\cal U}).$$
\end{theorem}
Our analysis significantly improves over the previous best known bound of $O(\sqrt m)$ for the performance of affine policies for budget of uncertainty sets. Furthermore, since our performance bound matches the hardness of approximation, affine policy provides an optimal approximation for~\eqref{eq:ar} for budget of uncertainty sets. In particular, there is no polynomial time algorithm whose worst-case approximation is better than affine policies by more than a constant factor. Note that the above statements relate to the worst-case performance. For particular instances, it may be possible to get better solutions than affine policies.

\subsection{Construction of our affine solution.}

Our analysis is based on constructing a {\em good} approximate affine solution that has a worst-case cost being $O \left( \frac{\log n}{\log \log n} \right)$ times the optimal cost $z_{\sf AR}({\cal U})$. In this section, we present the construction of our affine solution and consequently prove Theorem \ref{thm:budget}. Let us first introduce the following notations.

\vspace{2mm}
\noindent
{\bf Notations.}
We consider an optimal solution $ \mb x^*,  \mb y^*( \mb h) $ where $ \mb h \in {\cal U}$, for the adjustable problem \eqref{eq:ar}. Let $\opt$ be the optimal cost for \eqref{eq:ar} and ${\opt}_1, {\opt}_2$ respectively the first stage cost and the second stage cost associated with $ \mb x^*,  \mb y^*( \mb h) $, i.e. 
\begin{align*}
&{\opt}_1 = \mb c^T \mb x^*  \\
&{\opt}_2 = \max_{ \mb h \in {\cal U}} \; \mb d^T \mb y^* ( \mb h ) \\
 &{\opt}  ={\opt}_1+{\opt}_2 = z_{\sf AR}({\cal U}).
\end{align*}
We would like to remark that this split is not unique since there might be other optimal solutions for \eqref{eq:ar}. \red{ For all $i \in [m]$, we denote
$$ \alpha_i = 1 -( \mb A \mb x^*)_i$$
In particular, if $\alpha_i$ is negative the first stage solution $\mb x^*$ covers the full unit requirement in the $i$-th component, if not $\alpha_i$ corresponds  to the remaining requirement that needs to be covered eventually by the second-stage solution $\mb y^*(\cdot)$.}

We denote $z(\mb h )$, the cost of covering the requirement $\mb h$ in the second-stage, i.e.  
\begin{equation} \label{eq:z-cost}
 z(\mb h ) = \min_{\mb y \geq \mb 0} \left\{ \mb d^T \mb y \; \bigg\vert \; \mb B \mb y \geq \mb h \right\}.
\end{equation}
We refer to problem \eqref{eq:z-cost} as the \textit{fractional covering problem}. For  any ${\cal W} \subseteq [m]$, we denote $\mathbbm{1}({\cal W}) \in \mathbb{R}^m$ the indicator of ${\cal W}$, i.e. 
$$\mathbbm{1}({\cal W})= \sum_{i \in {\cal W}} \mb e_i.$$ 
For simplicity, we use the following notation 
$$z(\mathbbm{1}({\cal W}) )= z( {\cal W}).$$

\vspace{2mm}
\noindent
{\bf Our construction.}
For all $i \in [m]$, recall $z(\mb e_i)$  the optimal cost to cover component $\mb e_i$ in the second stage as defined in \eqref{eq:z-cost}. Let ${\mb v}_i$ be the optimal corresponding solution, i.e.,
$$ \mb v_i \in  \argmin_{\mb y \geq \mb 0} \left\{ \mb d^T \mb y \; \bigg\vert \;   \mb B \mb y \geq  \mb e_i  \right\}.$$
We split the components $\{ 1,2, \ldots,m\}$ based on a threshold into two sets ${\cal I}$ and its complement ${\cal I}^c$:
\begin{align*}
& {\cal I} = \left\{ i \in [m] \; \bigg\vert \;  \alpha_i  >0  \;\; \text{and}  \; \; \frac{ \alpha_i z(\mb e_i)}{w_i}   \leq \beta \cdot \opt \right\} \\
&{\cal I}^c =[m] \setminus  {\cal I},
\end{align*}
\red{where $$ \beta =   \frac{4 \log n}{\log \log n} .$$}
\red{We cover  a fraction of ${\cal I}$ (\textit{inexpensive components}) using a linear solution and  the remaining fraction of $\cal I$ along with  ${\cal I}^c$ (\textit{expensive components}) using a static solution. }

\vspace{2mm}
\noindent
{\bf Linear part.} We cover \red{a fraction of }the components of ${\cal I}$ using the following linear solution for any $ \mb h \in {\cal U}$,
\begin{equation} \label{eq:linear-sol}
 \red{\mb y_{\sf Lin}( \mb h )   = \sum_{i \in {\cal I} }  \alpha_i  h_i {\mb v}_i .}
\end{equation}

\vspace{2mm}
\noindent
{\bf Static part.}   \red{We use a static solution to cover the remaining components $\mb e_i$ where $ i \in {\cal I}^c$ and $(1-\alpha_i)^+ \mb e_i$ for $ i \in {\cal I}$ }. In particular, we consider the following static problem
\begin{equation} \label{eq:static}
\red{ (\mb x_{\sf Sta}, \mb y_{\sf Sta}) \in \argmin_{\mb x \in {\cal X}, \mb y \geq \mb 0} \left\{ \mb c^T \mb x+ \mb d^T \mb y \; \bigg\vert \;  \mb A \mb x+ \mb B \mb y \geq \sum_{ i \in {\cal I}^c} \mb e_i     + \sum_{ i \in {\cal I}}  (1-\alpha_i)^+ e_i \right\},}
\end{equation}
and  denote  $$z_{\sf Sta}=\mb c^T \mb x_{\sf Sta}+ \mb d^T \mb y_{\sf Sta}.$$ Therefore our candidate affine solution is
\begin{equation}\label{eq:feasible-affine}
\begin{aligned}
 \mb x  & =\mb x_{\sf Sta}  \\
 \mb y ( \mb h ) &   =  \mb y_{\sf Lin}( \mb h ) +  \mb  y_{\sf Sta}, \; \forall \mb h \in{\cal U}. \\
\end{aligned}
\end{equation}

\vspace{2mm}
\noindent
{\bf Feasibility.} We first show that  our candidate solution \eqref{eq:feasible-affine} is feasible for the adjustable problem \eqref{eq:ar}. The proof is a direct consequence of our construction. In particular, we have the following lemma.

\begin{lemma} \label{lem:feasibilty}
The affine solution in \eqref{eq:feasible-affine}  is feasible for the adjustable problem \eqref{eq:ar}.  
\end{lemma}
\red{\proof{Proof.}
We have, 
$$ \mb B  \mb y_{\sf Lin}( \mb h )  =  \sum_{i \in {\cal I}} \alpha_i h_i \mb B {\mb v}_i  \geq \sum_{i \in {\cal I}} \alpha_i h_i \mb e_i,$$ 
and
$$ \mb A \mb x_{\sf Sta} + \mb B \mb y_{\sf Sta}  \geq   \sum_{ i \in {\cal I}^c} \mb e_i     + \sum_{ i \in {\cal I}}  (1-\alpha_i)^+ e_i
\geq   \sum_{ i \in {\cal I}^c} h_i \mb e_i     + \sum_{ i \in {\cal I}}  h_i (1-\alpha_i)^+ e_i $$
 where the last inequality holds because $ h_i \in [0,1]$ for all $i \in [m]$.
Therefore, the solution in \eqref{eq:feasible-affine} verifies
\begin{align*}
&\mb A \mb x + \mb B \mb y( \mb h )  \geq \sum_{i \in {\cal I}^c} h_i \mb e_i+    \sum_{i \in {\cal I}}     ((1-\alpha_i)^+ + \alpha_i)               h_i \mb e_i  \geq \mb h. \\
  & \red{\mb{x}  \in  {\cal X} }\\
& \mb{y}(\mb{h}) \geq \mb 0, \;  \; \; \forall \mb h \in {\cal U}.
\end{align*}
 }
\hfill
\Halmos
\endproof

\red{ We would like to remark that the construction of the linear ans static parts  not only depends on the uncertainty set ${\cal U}$, but also depends
on all the parameters of the instance, i.e., $\mb A, \mb B, \mb c, \mb d$.}  This is in contrast to the analysis in \cite{BG10} where the construction of affine policies depends only on $\cal U$.

\subsection{Performance analysis.}
 We analyze separately the cost of the static and linear parts. For the linear part, the cost analysis is a direct consequence of our construction. In fact, we leave only  inexpensive scenarios to the linear part, \red{ i.e., scenarios $\alpha_i \mb e_i$ such that $\alpha_i z(\mb e_i)/w_i$ is less than the threshold $\beta \cdot \opt$.} We know that for all $\mb h \in {\cal U}$, we have
$\sum_{i=1}^m w_i h_i \leq 1$. Hence, the cost of linear part is at most $\beta \cdot \opt$. In particular, we have the following lemma.


\begin{lemma}[\bf {Cost of Linear part}] \label{lem:cost-linear}
The cost of the linear part $\mb y_{\sf Lin}(\mb h)$ defined in \eqref{eq:linear-sol}  is at most $\beta \cdot \opt $ for any $\mb h \in {\cal U}$.
\end{lemma}

\proof{Proof.}
\red{
We have for all $\mb h \in {\cal U}$,
$$
\mb d ^T\mb y_{\sf Lin}( \mb h )   =   \sum_{i \in {\cal I}} \alpha_i h_i \mb d^T{\mb v}_i =   \sum_{i \in {\cal I}} \alpha_i h_i  z( \mb e_i)   \leq   \beta \cdot \opt  \cdot \sum_{i \in {\cal I}} w_i h_i   \leq   \beta \cdot \opt,\\
$$
where the first inequality holds because $\alpha_i z( \mb e_i)  \leq w_i \beta \cdot \opt$ for all $i \in {\cal I}$ and the second inequality follows as $ \sum_{i \in {\cal I}} w_i h_i \leq 1$ for all $\mb h \in {\cal U}$.}
\hfill
\Halmos
\endproof

The key part is to analyze the cost of the static part. In fact, we show that the cost of the static part is also $O( \beta) \cdot \opt $. This relies on a structural result on  fractional covering problems. Intuitively, let us explain the structural result in the special case \eqref{eq:k-ones} of the budget of uncertainty set \red{and $\alpha_i$ are 0 or 1.}
We show that if the cost of covering every single component $\mb e_i, $ for $ i  \in {\cal J}$ is \textit{expensive}, i.e.  $ z(\mb e_i) >\beta \cdot \opt /k   $ and the cost of covering any $k$  components is \textit{inexpensive}, i.e. less than $2 \opt$. Then, the cost of covering all components of ${\cal J}$ is not too costly and can not exceed $\beta \cdot \opt$. The formal general statement is given in the following lemma.
\begin{lemma}[\bf{Structural Result}] \label{lem:structural-resulto}
Consider $\mb B \in \mathbb{R}^{m\times n}_+ $,  $ \mb d \in \mathbb{R}^n_+$ and  ${\cal J} \subseteq [m]$. Let $z(\mb h)$ be the cost of covering $\mb h$ as defined in~\eqref{eq:z-cost}.
Suppose there exists    $ \gamma >0$ and $w_i \in (0,1], \; \forall i \in {\cal J}$  such that the following two conditions are satisfied:
\begin{enumerate}
\item for all $ i \in {\cal J}$,  $$\frac{ z(\mb e_i)}{w_i}        >   4 \gamma \cdot \frac{\log n}{\log \log n},  $$ 
\item  for all ${\cal W} \; \subseteq {\cal J}$,
$$\sum_{i \in  {\cal W}} w_i \leq 1 \; \mbox{ implies }\; z({\cal W}) \leq \gamma.$$
\end{enumerate}
Then,  $$z({\cal J}) \leq  4 \gamma \cdot \frac{\log n}{\log \log n}.$$
\end{lemma}

\red{Lemma \ref{lem:structural-resulto} is a generalization of the result in Gupta et al. \cite{gupta2014thresholded}  for the set covering problem (see Theorem 7.1 in \cite{gupta2014thresholded}).  In particular, our result hold for any constraint matrix $ \mb B$ and a budgeted set with general  $w_i$ for $i \in [m]$, whereas the Gupta et al. [22], discuss the special case  when $\mb B \in \{0,1\}^{m \times n}$,  and a budget of uncertainty set with $w_i = 1/k$. Furthermore, we improve the approximation bound from $O(\log n)$ in \cite{gupta2014thresholded} to $O(\log n/\log \log n)$.} We present the proof of Lemma \ref{lem:structural-resulto} later in Section~\ref{section:claims}. But let us first use the structural result to show that the cost of the static part is $O(\beta ) \cdot \opt$ and consequently prove Theorem \ref{thm:budget}.  In particular, we have the following lemma.

\begin{lemma}[\bf{Cost of Static part}] \label{lem:cost-static}
The cost $z_{\sf Sta} $ of the static part $(\mb x_{\sf Sta}, \mb y_{\sf Sta})$ defined  in \eqref{eq:static} is  $O (\beta) \cdot \opt$. 
\end{lemma}

\red{
\proof{Proof.}
Consider the following sets
$${\cal J}_{1}= \left\{ i \in [m] \; \bigg\vert \; \alpha_i \leq  0 \right\}.$$
and 
$$ {\cal J}_2= {\cal I}^c  \setminus{\cal J}_1 =\left\{ i \in [m] \; \bigg\vert \;  \alpha_i  >0  \;\; \text{and}  \; \; \frac{ \alpha_i z(\mb e_i)}{w_i}   >\beta \cdot \opt \right\}.$$
For $i=1,\ldots,m$,  denote $\mb B^T_i$ the $i$-th row of $\mb B$ and let $  \mb{\tilde{B}}^T_i = \mb B^T_i / \alpha_i$. We have for $i \in [m]$, 
$$ \alpha_i z( \mb e_i)= \min_{\mb y \geq \mb 0} \left\{ \mb d^T \mb y \; \bigg\vert \; \mb{\tilde{B}} \mb y \geq \mb e_i \right\}.$$
We apply the structural Lemma \ref{lem:structural-resulto} with the parameters $\mb{\tilde{B}}, \mb d , {\cal J}_{2}$ and $\gamma =  \opt$. Let us verify the assumptions of Lemma~\ref{lem:structural-resulto}. For all $i \in {\cal J}_{2}$, we have
$$  \frac{  \alpha_i  z( \mb e_i)}{w_i}  > \beta  \cdot \opt = 4 \gamma \cdot \frac{\log n}{\log \log n}.$$
For any ${\cal W} \subseteq {\cal J}_{2}$  such that $\sum_{i \in {\cal W}} w_i \leq 1$, we have  $\mb h = \mathbbm{1} ({\cal W} ) \in \cal{U}$. By feasibility of the optimal solution, we know that 
$$ \mb A \mb x^*+ \mb B \mb y^*(\mb h) \geq \mb h,$$
which implies,
$$ \mb B \mb y^*(\mb h) \geq \sum_{i \in {\cal W}}  \mb e_i -       \sum_{i=1}^m (1- \alpha_i) \mb e_i  = \sum_{i \in {\cal W}} \alpha_i \mb e_i + \sum_{i \notin {\cal W}} ( \alpha_i-1) \mb e_i   $$
In particular, we have for all $i \in {\cal W}$, $(\mb B \mb y^*(\mb h))_i \geq \alpha_i$. Moreover, $\mb B$ and $\mb y^*$ are non-negative which implies that
$$ \mb B \mb y^*(\mb h) \geq \sum_{i \in {\cal W}} \alpha_i \mb e_i,$$
and therefore,
$$ \tilde{\mb B}\mb y^*(\mb h) \geq \sum_{i \in {\cal W}} \mb e_i = \mb h.$$
This means that $\mb{y}^* (\mb h)$ is a feasible solution for the covering problem \eqref{eq:z-cost} with constraint matrix $\tilde{\mb B}$ and requirement $\mb h$. Therefore,
$$\min_{\mb y \geq \mb 0} \left\{ \mb d^T \mb y \; \bigg\vert \; \mb{\tilde{B}} \mb y \geq \mb h \right\} \leq  \mb d^T  \mb{y}^* (\mb h)  \leq   \opt_2 \leq    \opt = \gamma.$$
This verifies the second assumption of Lemma \ref{lem:structural-resulto}. Therefore, from the structural result, we have 
$$ \min_{\mb y \geq \mb 0} \left\{ \mb d^T \mb y \; \bigg\vert \; \mb{\tilde{B}} \mb y \geq \sum_{i \in {\cal J}_{2} }  \mb e_i\right\}   \leq 4 \gamma  \frac{\log n}{\log \log n} =\beta \cdot \opt.$$ 
Denote $\mb y_2$ an optimal solution corresponding to the above minimization problem. In particular, we have $\mb d^T \mb y_2 \leq  \beta \cdot \opt$ and 
$$ \mb B \mb y_2 \geq \sum_{i \in {\cal J}_{2} } \alpha_i \mb e_i.$$
Furthermore, by feasibility of the optimal solution for \eqref{eq:ar}, we have
$$ \mb A \mb x^*+ \mb B \mb y^*(\mb 0) \geq \mb 0,$$
and  we know that $\mb B \mb y^*(\mb 0) \geq \mb 0$. This  implies,
$$ \mb B \mb y^*(\mb 0) \geq \sum_{ i=1}^m (\alpha_i -1)^+ \mb e_i.$$
Putting all together, we have
\begin{align*}
\mb A \mb x^* + \mb B\mb y^*(\mb 0)+ \mb B \mb y_2 &\geq \sum_{i=1}^m (1-\alpha_i) \mb e_i + \sum_{ i=1}^m (\alpha_i -1)^+ \mb e_i + \sum_{i \in {\cal J}_{2} } \alpha_i \mb e_i \\
&= \sum_{i=1}^m (1-\alpha_i)^+ \mb e_i  + \sum_{i \in {\cal J}_{2} } \alpha_i \mb e_i \\
&= \sum_{i \in {\cal I}}    (1-\alpha_i)^+ \mb e_i  +   \sum_{i \in {\cal J_1}}     (1-\alpha_i)^+ \mb e_i+   \sum_{i \in {\cal J_2}}    (1-\alpha_i)^+\mb e_i        +                 \sum_{i \in {\cal J}_{2} } \alpha_i \mb e_i \\
&\geq  \sum_{i \in {\cal I}}    (1-\alpha_i)^+ \mb e_i  +   \sum_{i \in {\cal J_1}}      \mb e_i+   \sum_{i \in {\cal J_2}}    \mb e_i \\
& =  \sum_{i \in {\cal I}}    (1-\alpha_i)^+ \mb e_i  +   \sum_{i \in {\cal I}^c}      \mb e_i
\end{align*}
where the last inequality holds because $ \alpha_i \leq 0 $ for all $ i \in {\cal J}_1$ and $   (1-\alpha_i)^+ + \alpha_i \geq 1$ for all $ i \in  {\cal J}_2$. Moreover, $\mb x^* \in {\cal X}$ and $(\mb y^*(\mb 0)+  \mb y_2)$ is non-negative.  Hence, we have $(\mb x^*, \mb y^*(\mb 0)+  \mb y_2         )$ is a feasible solution for the static problem in \eqref{eq:static}, therefore
$$ z_{\sf Sta} \leq     \mb c^T \mb x^* + \mb d^T  y^*(\mb 0)  + \mb d^T     \mb{y}_2  
 \leq   \opt+   \beta \cdot  \opt = O( \beta) \cdot \opt.$$
\hfill
\Halmos
\endproof
}

\proof{Proof of Theorem \ref{thm:budget}.}
Lemma \ref{lem:feasibilty} show that our affine solution \eqref{eq:feasible-affine} is feasible for the adjustable problem \eqref{eq:ar}. Lemma \ref{lem:cost-linear} and \ref{lem:cost-static} show that  the cost of the affine solution \eqref{eq:feasible-affine}  is less than $\beta  \cdot\opt +O( \beta) \cdot \opt = O(\beta ) \cdot \opt $ which implies that $$ z_{\sf Aff}({\cal U}) = O \left( \frac{\log n}{\log \log n} \right)\cdot z_{\sf AR}({\cal U}).$$
\hfill
\Halmos
\endproof

\subsection{Proof of the Structural Result.} \label{section:claims}

\red{We give a proof by contradiction. 
The assumptions in Lemma~\ref{lem:structural-resulto} can be interpreted as follows. Let $\eta  =  \frac{4 \log n }{\log \log n}$. The first assumption states that the cost of covering any single component, $\mb e_i$ for $i \in {\cal J}$ is large (at least $w_i \eta \cdot \gamma$). The second assumption states that the cost of any feasible (integral) scenario, ${\cal W} \subseteq {\cal J}$ with $\sum_{i \in {\cal W}} w_i \leq 1$ is at most $\gamma$. We need to show that the cost of covering all the components, ${\cal J}$ is at most $\eta\cdot \gamma$. For the sake of contradiction, suppose that  $z({\cal J})  > \eta \gamma$. We will construct a feasible scenario, ${\cal W} \subseteq {\cal J}$ (where $h_i =1$ for all $i \in {\cal W}$ and $\sum_{i \in {\cal W}} w_i \leq 1$) where the cost, $z({\cal W}) > \gamma$ violating the second assumption in Lemma~\ref{lem:structural-resulto}. To construct this scenario, we consider the dual of the primal covering problem $z({\cal J})$. The dual is a packing problem where the ratio of the right hand sides and the constraint coefficients is {\em large} (from the first assumption in the lemma). This allows us to construct an approximate {\em integral} dual solution of the problem (using randomized dual rounding) where we lose only a factor $\eta$ in the objective value as compared to the optimal (fractional) dual solution. We then use this approximate integral dual to construct a scenario ${\cal W}$ with cost greater than $\gamma$ that gives us the contradiction. Details of the proof are provided below.}

For all $k \in  {\cal J} $, recall ${\mb v}_k$ the optimal solution corresponding to $z(\mb e_k)$. We have  $\|  {\mb v}_k \|_0 =1$, i.e. $ z(\mb e_k) = d_{\ell} { v}_{k\ell}$ where 
$$ \ell = \underset{ \underset {B_{kj} \neq 0} {j=1,\ldots,n}} \argmin  \; \frac{d_j}{B_{kj}}.$$
In particular, we have for all $j \in [n]$ such that $B_{kj} \neq 0$, 
\begin{align*}
\frac{d_j}{B_{kj}} \geq \frac{d_{\ell}}{B_{k{\ell}}} = d_{\ell} v_{k\ell} = z(\mb e_k) > \eta \gamma w_k,
\end{align*}
i.e., for all $j \in [n]$,
\red{$$d_j \geq \eta \gamma \cdot \max_{k \in {\cal J}}  (w_k B_{kj} )  .$$} 
For $j \in [n]$, denote
\red{$$\hat{d}_j = \frac{d_j}{ \eta \gamma \cdot \max_{k \in {\cal J}} (w_k B_{kj})}, $$}
and  for all $i \in {\cal J} ,j \in [n]$,
\red{$$ \hat{B}_{ij}= \frac{w_i B_{ij}}{\max_{k \in {\cal J}} (w_k B_{kj})}.$$}
In particular, we have for all $i \in {\cal J} ,j \in [n]$, $\hat{B}_{ij} \in [0,1]$ and for all $j \in [n]$, $\hat{d}_j >1$. For any $ {\cal W} \subseteq {\cal J}$, consider the following problem
\begin{equation}\label{eq:one-stage-hat}
 \hat z({\cal W})  = \min_{\mb y \geq \mb 0 } \left\{ \mb{\hat{d}}^T \mb y \; \bigg\vert \; \hat{\mb B} \mb y \geq \sum_{ i \in {\cal W}}  w_i \mb  e_i   \right\}.
\end{equation}
We show that for any $ {\cal W} \subseteq {\cal J}$, $ \hat{z}({\cal W} )$ is just a scaling of $z({\cal W}) $. In particular, we have, 

\begin{claim}\label{claim:3.1} 
$z({\cal W}) = \eta \gamma \cdot  \hat{z}({\cal W} ) .$
\end{claim}

\vspace{1mm}
\noindent
We present the proof of Claim~\ref{claim:3.1} in Appendix~\ref{appendix:lem:structural-resulto}. To show that $z({\cal J}) \leq \eta \gamma$\red{, we} show equivalently that   $  \hat{z}({\cal J} )  \leq 1 $. For the sake of contradiction, suppose that $ \hat{z}({\cal J})  >1$. Our goal is to construct a scenario that contradicts condition 2 of the lemma. We use ideas on dual rounding and randomized solutions from \cite{FJMM07} and \cite{gupta2014thresholded}.
In particular, let the dual  problem of $\hat{z}({\cal J})$ be

\begin{equation}\label{eq:dual}
\hat{\Delta}_{\cal J}= \max_{\mb z \geq \mb 0}  \left\{ \sum_{ i \in {\cal J}} w_i z_i \; \bigg\vert \;   \sum_{ i \in {\cal J}} \hat{B}_{ij} z_i \leq \hat{d}_j \; \; \forall j \in [n]   \right\}.
\end{equation}
Denote $ \mb{z^*}$ the optimal solution for the dual problem \eqref{eq:dual}. By strong LP duality, we have $$ \hat{\Delta}_{\cal J} = \hat{z}({\cal J}),$$ and therefore,    
$$ \hat{\Delta}_{\cal J}=  \sum_{i \in {\cal J}}  w_i z^*_i > 1.$$
We define the following randomized solution for all $ i \in {\cal J}$,
$$ Z_i = \lfloor{z_i^*}\rfloor + {\xi}_i,  $$
where ${\xi}_i,$ for $ i \in {\cal J}$ are independent Bernoulli variables with parameter $z_i^*- \lfloor{z_i^*}\rfloor$, i.e., 
$$ {\xi}_i ={\sf  Ber} ( z_i^*- \lfloor{z_i^*}\rfloor ).$$

\begin{claim}\label{claim:3.2} 
With probability at least $1- O(1/n)$, $$\left( \frac{ 2 Z_i}{\eta},{ i \in {\cal J}}\right),$$ is a feasible solution to the dual problem~\eqref{eq:dual}.
\end{claim}
\vspace{2mm}
\noindent 
We show that  $  \left( \frac{ 2 Z_i}{\eta},{ i \in {\cal J}}\right)$ satisfies the constraints of \eqref{eq:dual} with high probability by using Chernoff bound concentration inequalities. The proof of Claim~\ref{claim:3.2} is presented in Appendix \ref{appendix:lem:structural-resulto}. Furthermore, we show that the cost of our randomized solution   $\left( \frac{ 2 Z_i}{\eta},{ i \in {\cal J}}\right)$ is greater than $\frac{1}{\eta}$ with a constant probability. In particular, we have the following claim.

\begin{claim}\label{claim:3.3} 
$   \mathbb{P} \left(       \sum_{i \in {\cal J}}w_i Z_i  > \frac{1}{2} \right) \geq 1 - e^{  -\frac{ 1}{8}}$.
\end{claim}

We use a concentration bound to prove Claim~\ref{claim:3.3}. The proof is presented in Appendix \ref{appendix:lem:structural-resulto}. Putting Claim~\ref{claim:3.2} and Claim~\ref{claim:3.3} together,  we have that 
$$  \left( \frac{ 2 Z_i}{\eta},{ i \in {\cal J}}\right),$$ 
is feasible for \eqref{eq:dual}  with high probability and has a cost $\sum_{i \in {\cal J}}w_i \frac{ 2 Z_i}{\eta} $ strictly greater than $ \frac{1}{\eta}$ with a non-zero constant probability. Therefore, there exists a deterministic solution for problem \eqref{eq:dual} with a cost at least  $\frac{1}{\eta}$. For simplicity of notations, let us assume that $  \left( \frac{ 2 Z_i}{\eta},{ i \in {\cal J}}\right)$ is such a solution. Let us order $w_i Z_i$ in an increasing order, i.e.,
$$ w_{(1)} Z_{(1)} \geq w_{(2)} Z_{(2)} \geq \ldots \geq w_{(\vert {\cal J} \vert )} Z_{(\vert {\cal J} \vert)}.$$
We know that $\sum_{i \in {\cal J}} w_i Z_i > \frac{1}{2}$. Denote $L$ the index such that 
$$ \sum_{i =1}^{L-1}w_{(i)} Z_{(i)}   \leq \frac{1}{2} \qquad \text{and} \qquad \sum_{i =1}^{L}w_{(i)}  Z_{(i)}  > \frac{1}{2}.$$
Note that $Z_{(i)}  $ are integral and $Z_{(L)} \neq 0$. Hence for all $i=1, \ldots, L$, $Z_{(i)} \geq 1$. Therefore,
$$ \sum_{i =1}^{L-1}w_{(i)}  \leq  \sum_{i =1}^{L-1}w_{(i)}Z_{(i)}  <  \frac{1}{2}.$$
Note that if $L=1$, $\sum_{i =1}^{L}w_{(i)}   = w_{(1)} \leq 1$ because all $w_i$ are in $[0,1]$. On the other hand, if $L \geq 2$, then 
$$w_{(L)} \leq w_{(L)}  Z_{(L)} \leq w_{(1)}  Z_{(1)}  \leq \sum_{i =1}^{L-1}w_{(i)} Z_{(i)}  < \frac{1}{2},$$
and therefore,
$$ \sum_{i =1}^{L}w_{(i)} = \sum_{i =1}^{L-1}w_{(i)}+ w_{(L)}   \leq \frac{1}{2}+ \frac{1}{2 }=1.$$
Therefore, in both cases we have $\sum_{i =1}^{L}w_{(i)}  \leq 1.$
Denote ${\cal W} \subseteq {\cal J} $ the set of indices corresponding  to the top $L$ of $w_{(i)} Z_{(i)}$. In particular, we have, 
$$ \sum_{i \in {\cal W}} w_i \leq 1. $$
Note that $ \sum_{i \in {\cal W}} w_i Z_i > \frac{1}{2},$ and consequently, $$ \sum_{i \in {\cal W}} w_i \frac{2 Z_i}{\eta} > \frac{1}{\eta}.$$  Consider the following problem 

\begin{equation}\label{eq:dualW}
\hat{\Delta}_{\cal W}= \max_{\mb z \geq \mb 0}  \left\{ \sum_{ i \in {\cal W}} w_i z_i \; \bigg\vert \;   \sum_{ i \in {\cal W}} \hat{B}_{ij} z_i \leq \hat{d}_j \; \; \forall j \in [n]   \right\}
\end{equation}
and its dual,
\begin{equation}\label{eq:prob}
 \hat z({\cal W})  = \min_{\mb y \geq \mb 0} \left\{ \mb{\hat{d}}^T \mb y \; \bigg\vert \; \hat{\mb B} \mb y \geq \sum_{ i \in {\cal W}} w_i \mb e_i \right\}.
\end{equation}
 We have shown the existence of a  solution $(\frac{2Z_i}{\eta})_{i\in {\cal W}}$ to problem \eqref{eq:dualW} with a cost strictly greater than $\frac{1}{\eta}$. Hence, by LP duality 
$$  \hat z({\cal W})  >\frac{1}{\eta}.$$ 
Note that $  z({\cal W})  =  \eta \gamma\cdot \hat z({\cal W}) $. Hence, $    z({\cal W})  >  \gamma$ which contradicts  condition 2 of our lemma.

\section{Intersection of Disjoint budget constraints.}
In this section, we consider more general uncertainty sets that are defined by intersection of budget constraints that model many practical settings. We first consider the case where the budget constraints are {\em disjoint}. In particular, consider $S_1, S_2, \ldots, S_L$ a partition of $\{1,2,\ldots,m\}$, i.e. 
$$\bigcup_{{\ell}=1}^L S_{\ell} = [m] \qquad  \text{and} \qquad  S_i \cap S_j = \emptyset, \; \forall i \neq j.$$ We define  {\em disjoint constrained budgeted sets} as follows
\begin{equation}\label{eq:partition-matroid}
 {\cal U} = \left\{ \mb h \in [0,1]^m  \; \bigg\vert \;   \sum_{i \in S_{\ell}} w_{\ell i} h_i \leq 1  \;\;  \forall \ell \in [L] \right\}.
\end{equation}
where $S_{\ell}, \ell=1, \ldots,L$ is a partition of $[m]$.
This is an important class of uncertainty sets that generalizes the budget of uncertainty set \eqref{eq:budget}. These are essentially Cartesian product of $L$ budget of uncertainty sets. A special case of this class of uncertainty sets  where all $w_{\ell i}$ are equal has been considered for example in Gupta et al. \cite{gupta2016robust} and Feige et al. \cite{FJMM07}. Recall for $L=1$, namely the budget of uncertainty set \eqref{eq:budget}, affine policy gives the optimal approximation to \eqref{eq:ar} (see Theorem \ref{thm:budget}). Our result in this section show that the performance of affine policy remains near-optimal for the more general class \eqref{eq:partition-matroid}.  In particular, we have the following theorem. 

\begin{theorem}\label{thm:partition-matroid}
Consider the two-stage adjustable problem \eqref{eq:ar} where $ {\cal U}$ is the disjoint constrained budgeted  set \eqref{eq:partition-matroid}.  Then,
$$ z_{\sf Aff}({\cal U}) = O\left( \frac{\log^2 n }{\log \log n} \right)\cdot z_{\sf AR}({\cal U}) .$$
\end{theorem}

Our analysis relies on constructing a feasible affine solution for \eqref{eq:ar} and relating the worst-case performance to a lower bound of \eqref{eq:ar}. In particular, we consider the {\em online fractional covering problem} and use an online algorithm with $O(\log n)$-{\em competitive ratio} to both construct a feasible affine and also a lower bound for \eqref{eq:ar}. The performance bound of our feasible affine is related to the competitive ratio of the online algorithm. We first introduce some preliminaries before discussing our construction and analysis.


\subsection{Online fractional covering.}
\red{Recall, for $i=1,\ldots,m$,   $\alpha_i = 1 - ( \mb A \mb x^*)_i$. We consider $\tilde{\mb B}  \in \mathbb{R}_+^{m \times n}$  where the $i$-th row $  \mb{\tilde{B}}^T_i = \mb B^T_i / \alpha_i$ if $\alpha_i > 0$ and $  \mb{\tilde{B}}^T_i = \mb B^T_i $ otherwise. Note that $\tilde{\mb B}$ is a non-negative matrix. 
We consider the (offline) fractional covering problem 
$$ \Theta( \mb h) = \min_{\mb y \geq \mb 0} \left\{ \mb d^T \mb y \; \big\vert \; \tilde{ \mb B} \mb y \geq \mb h \right\},$$
for any requirement $\mb h \in \{0,1\}^m$. The \textit{online fractional covering problem} is an online version of the covering problem where the requirements are revealed online in a sequential manner. In particular, at each step we get a new constraint $ \sum_{j=1}^n \tilde{B}_{ij} y_j \geq 1$ for some $i$ and the algorithm needs to augment the current solution to satisfy the new requirement in each step.
}
This problem has been studied in the literature.  We refer the reader to Buchbinder and Naor \cite{buchbinder2009design} for an extensive discussion of the problem. Buchbinder and Naor \cite{buchbinder2009online} give an online algorithm $\cal A$ for the online fractional covering problem that is $O(\log n)$-\textit{competitive} (see Theorem 4.1 in \cite{buchbinder2009online}). In other words, the cost of the solution given by $\cal A$ for any set and sequence of requirements is at most $O(\log n)$ times the  cost of the optimal solution of the corresponding offline covering problem where all the requirements are known upfront. In particular, for any sequence of requirements $\tau$, we have
$$ \max_{\tau} \; \frac{{\cal A}(\tau)}{\Theta(\tau )} =O( \log n),$$ 
where ${\cal A}(\tau)$ is the online covering cost and $\Theta(\tau)$ is the offline covering cost. Note that the competitive ratio guarantee also holds for the case where in each step, we get a subset of constraints instead of just a single constraint. We consider the algorithm $\cal A$ of Buchbinder and Naor~\cite{buchbinder2009online} for our analysis. Let us introduce some notations that we will use for our construction and analysis. 

\vspace{2mm}
\noindent
{\bf Notations.}
Consider a sequence of subsets of constraints given by $(S_1,\ldots,S_R)$ where $S_r \subseteq [m]$ and $S_r \cap S_{r'} = \emptyset$ for all $r\neq r'$. In particular, in step $r$ we get subset $S_r$ of constraints
\red{$$  \sum_{j=1}^n \tilde{B}_{ij} y_j \geq 1 \qquad \forall i \in S_r.$$}
For brevity of notations, let $$\mb h_r = \mathbbm{1}(S_r).$$
In particular, the sequence $( \mb h_1, \mb h_2,\ldots, \mb h_R)$ verifies $\mb h_r\in \{0,1\}^m$ for all $ r \in [R]$ and  $\sum_{r=1}^R \mb h_r \leq \mb e.$
 We introduce the following definitions.

\vspace{2mm}
\begin{enumerate}
\item {\bf Online cost.} We denote ${\cal A} ( \mb h_1, \mb h_2,\ldots, \mb h_r)$ the  (online) cost of covering the sequence $( \mb h_1, \mb h_2,\ldots, \mb h_r)$ using the online algorithm ${\cal A}$.

\vspace{2mm}
\item {\bf Online augmenting cost.} We denote ${\cal A} ( {\mb h}_{r+1}  \vert \; \mb h_1, \mb h_2,\ldots, \mb h_r)$ the extra cost to cover ${\mb h}_{r+1} $ using the online algorithm $\cal A$ when the algorithm have already covered  the sequence $( \mb h_1, \mb h_2,\ldots, \mb h_r)$. By definition, the online augmenting cost is given by
\begin{equation} \label{eq:aug}
{\cal A} ( {\mb h}_{r+1} \vert \mb h_1, \mb h_2,\ldots, \mb h_r) ={\cal A} (  \mb h_1, \mb h_2,\ldots, \mb h_r , {\mb h}_{r+1} )-  {\cal A} (\mb h_1, \mb h_2,\ldots, \mb h_r).
\end{equation}

\vspace{2mm}
\item {\bf Greedy augmenting cost.} We denote ${\sf Aug} ( \mb h_{r+1} \vert \; \mb h_1, \mb h_2,\ldots, \mb h_r)$ the optimal cost to cover $\mb h_{r+1}$  given that the sequence $( \mb h_1, \mb h_2,\ldots, \mb h_r)$ was already covered by the online algorithm $\cal A$. In particular, the greedy augmenting cost is given by
\red{\begin{equation} \label{eq:optaug}
{\sf Aug} ( \mb h_{r+1} \vert \; \mb h_1, \mb h_2,\ldots, \mb h_r) = \min_{ \mb y \geq \mb 0} \left\{ \mb d^T \mb y \; \big\vert \; \tilde{ \mb B } \left( \mb y+ \mb y_r^{\cal A} \right)\geq \mb h_{r+1} \right\},
\end{equation}}
where $\mb y_r^{\cal A}$ is the online solution corresponding to the cost ${\cal A}( \mb h_1, \mb h_2,\ldots, \mb h_r)$.

\vspace{2mm}
\item  {\bf Offline cost.} \red{Denote $\Theta (\mb h_1, \mb h_2,\ldots, \mb h_r)$ the optimal (offline) cost to cover $( \mb h_1, \mb h_2,\ldots, \mb h_r)$ i.e.,
$$ \Theta(\mb h_1, \mb h_2,\ldots, \mb h_r) = \Theta \left( \sum_{ i =1}^r \mb h_i  \right) = \min_{\mb y \geq \mb 0} \left\{ \mb d^T \mb y \; \bigg\vert \; \tilde{ \mb B} \mb y \geq \sum_{ i =1}^r \mb h_i \right\}. $$}
\end{enumerate}
Since $\cal A$ is $O(\log n)$-competitive. Then for any sequence $( \mb h_1, \mb h_2,\ldots, \mb h_r)$,
\begin{equation} \label{thm:noga}
\red{ {\cal A} ( \mb h_1, \mb h_2,\ldots, \mb h_r)=O(\log n) \cdot   \Theta  \left( \sum_{ i =1}^r \mb h_i  \right).}
\end{equation}

\subsection{Construction of our affine solution.}

Similar to the proof of Theorem \ref{thm:budget}, we construct a feasible affine solution where we split the components of $[m]$ into two subsets and cover one using a linear solution and the remaining components  using a static solution. Consider $1,2,\ldots,L$ the blocks of components of the the disjoint constrained budgeted set \eqref{eq:partition-matroid}. For each block of components, we construct a threshold using the online fractional covering  algorithm $\cal A$. This threshold defines the  {\em expensive } components that we cover using a static solution and {\em inexpensive} components that we cover using a linear solution.

\vspace{2mm}
\noindent
{\bf Construction of  thresholds.} Let
\red{$$ 
\begin{aligned}
{\cal T} & = \left\{ i \in [m] \; \bigg\vert \; \alpha_i \leq 0 \right\}, \\
{\cal T}^c & =[m] \setminus {\cal T}.
\end{aligned}
$$}
 For $\ell =1 ,\ldots ,L$, denote $$\hat{S}_{\ell}= S_{\ell} \cap {\cal T}^c.$$ Let us define the following sets for all $\ell \in [L]$,
$$ {\cal U}_{\ell}=  \left\{ \mb h \in \{0,1\}^m  \; \Big\vert \;   \sum_{i \in \hat{S}_{\ell}} w_{\ell i} h_i \leq 1  \; \; \text{and} \; \;  h_i=0  \; \; \forall i  \notin \hat{S}_{\ell}  \right\}.$$
Note that $ \oplus_{\ell=1}^L {\cal U}_{\ell} \subseteq {\cal U}$. We construct a greedy sequence $(\mb a_1 , \mb a_2, \ldots , \mb a_L)$ where each $\mb a_{\ell}$ is chosen from some set ${\cal U}_{\ell}$ such that it maximizes the online augmenting cost \eqref{eq:aug} of the sequence.  Algorithm \ref{algo1} describes the procedure in details.
\begin{algorithm*}[]
\caption{Computing a greedy scenario $\mb a$}\label{algo1}
\begin{algorithmic}[1]
\State Initialize $ {\cal L}= \{ 1,2, \ldots, L \}$.
 \For{ ${\ell}=1,\ldots,L$}
 \State $(s,\mb b) = \underset{  s \in {\cal L}, \; \;  \mb b \in {\cal U}_s }{\argmax} \; \;  {\cal A}( \mb b \; \vert  \; \mb a_{1}, \mb a_{2}, \ldots , \mb a_{{\ell}-1} )    $
 \State Set $\mb a_{\ell}= \mb b $, update ${\cal L}= {\cal L}\setminus s$
 \EndFor
\end{algorithmic}
\end{algorithm*}
\red{
Algorithm~\ref{algo1} constructs the greedy sequence when the covering constraint matrix is $\mb{\tilde B}$ and therefore, the guarantee of the online algorithm ${\cal A}$ gives us a bound of $O(\log n)$ between the online cost, ${\cal A}(\mb a)$ and offline cost, $\Theta(\mb a)$ of the covering problem with $\tilde{\mb B}$ for the greedy sequence, $\mb a$. The following lemma relates this cost with {\sf OPT} for~\eqref{eq:ar}. 
\begin{lemma}\label{lem:lastone}
For all $ \ell \in [L]$, denote $\nu_{\ell}$ the cost of covering the sequence $(\mb a_1, \mb a_2, \ldots , \mb a_{\ell})$ using the online algorithm ${{\cal A}}$, i.e.,
$${\nu}_{\ell}= {\cal A}( \mb a_1, \mb a_2,\ldots, \mb a_{\ell}).$$ We have
$$ {\nu}_L = O(\log n) \cdot  \opt.$$
\end{lemma}
\proof{Proof.}
We suppose wlog that the sequence $(\mb a_1, \mb a_2, \ldots , \mb a_L)$ belongs respectively to  ${\cal U}_1,\; {\cal U}_2,\ldots,{\cal U}_L$ and let
$$ \mb a = \sum_{\ell=1}^L \mb a_{\ell} \in {\cal U}.$$
Then,
$$ \mb A \mb x^* + \mb B \mb y^*( \mb a) \geq \mb a.$$
Denote $\mb a = \sum_{i \in {\cal S}} \mb e_i $. Hence,
$$  \mb B \mb y^*( \mb a) \geq \sum_{i \in {\cal S}} \mb e_i     -  \sum_{i =1}^m  (1- \alpha_i) \mb e_i,        $$
i.e.,
$$  \mb B \mb y^*( \mb a) \geq \sum_{i \in {\cal S}}  \alpha_i \mb e_i     -  \sum_{i \notin {\cal S}}  (1- \alpha_i) \mb e_i.$$
Moreover, $\mb B$ and $ \mb y^*(\mb a)$ are non-negative. Hence, $\mb B \mb y^*( \mb a) \geq \mb 0$ and therefore,
$$  \mb B \mb y^*( \mb a) \geq \sum_{i \in {\cal S}}  \alpha_i \mb e_i,        $$
which is equivalent to
$$  \tilde{\mb B} \mb y^*( \mb a) \geq \sum_{i \in {\cal S}}   \mb e_i   = \mb a.$$
Therefore,
$$ \Theta(\mb a) =    \min_{\mb y \geq \mb 0} \left\{ \mb d^T \mb y \; \bigg\vert \; \mb{\tilde{B}} \mb y \geq \mb a \right\} \leq    \mb d^T y^*(\mb a)  \leq  \opt.$$
Finally,
$$ {\nu}_L = O(\log n) \cdot  {\Theta}(\mb a) = O(\log n) \cdot \opt,$$
where the first equality follows from \eqref{thm:noga}.
\hfill
\Halmos
\endproof
}

Now, we are ready to construct our feasible affine solution that has a cost  $O\left( \frac{\log^2 n}{\log \log n} \right)$ times $z_{\sf AR}( \cal U)$ using a linear and a static part.
Recall for all $i \in [m]$, $z( \mb e_i)$ the cost of covering component $\mb e_i$ in the second stage  as defined in \eqref{eq:z-cost} and $ {\mb v}_i $ an optimal corresponding solution. For all $\ell =1, \ldots,L$,  we consider the following sets of components
\red{$$ {\cal I}_{\ell}= \left\{ i \in S_{\ell} \; \bigg\vert  \;  \alpha_i >0 \; \; \text{and} \; \; \frac{\alpha_i z(\mb e_i)}{w_{\ell i}}  \leq \beta \cdot ({\nu}_{\ell}-{\nu}_{\ell-1} )\right\},$$}
where $$\beta=   \frac{8 \log n}{\log \log n}, $$ and $${\nu}_{\ell}-{\nu}_{\ell-1}= {\cal A} (  \mb a_{\ell} \vert \mb a_1, \mb a_2,\ldots, \mb a_{\ell-1}),$$ as defined in \eqref{eq:aug}. Denote
$${\cal I}= \bigcup_{\ell=1}^L {\cal I}_{\ell}.$$
and ${\cal I}^c$ its complement, i.e., ${\cal I}^c= [m] \setminus  {\cal I}$.

\vspace{2mm}
\noindent
{\bf Linear part.}
We cover \red{ a fraction of the components of ${\cal I}$} using the following linear solution for any $ \mb h \in {\cal U}$,
\begin{equation} \label{eq:linear-sol-matroid}
\red{ \mb y_{\sf Lin}( \mb h )   = \sum_{i \in {\cal I}}   \alpha_i h_i {\mb v}_i .}
\end{equation}

\vspace{2mm}
\noindent
{\bf Static part.}   \red{We use a static solution to cover the remaining components $\mb e_i$ where $ i \in {\cal I}^c$ and $(1-\alpha_i)^+ \mb e_i$ for $ i \in {\cal I}$ }. In particular, similar to \eqref{eq:static}  we consider the following static problem
\begin{equation} \label{eq:static-matroid}
\red{ (\mb x_{\sf Sta}, \mb y_{\sf Sta}) \in \argmin_{\mb x \in {\cal X}, \mb y \geq \mb 0} \left\{ \mb c^T \mb x+ \mb d^T \mb y \; \bigg\vert \;  \mb A \mb x+ \mb B \mb y \geq \sum_{ i \in {\cal I}^c} \mb e_i     + \sum_{ i \in {\cal I}}  (1-\alpha_i)^+ e_i \right\},}
\end{equation}
and  denote  $z_{\sf Sta} =\mb c^T \mb x_{\sf Sta} + \mb d^T \mb y_{\sf Sta}$. Our affine solution is given by
\begin{equation}\label{eq:feasible-affine-matroid}
\begin{aligned}
\mb x &  =\mb x_{\sf Sta}  \\
\mb y ( \mb h ) &  =  \mb y_{\sf Lin}( \mb h ) +  \mb  y_{\sf Sta}, \qquad \forall \mb h \in{\cal U}. \\
\end{aligned}
\end{equation}
We can show that the  affine  solution  \eqref{eq:feasible-affine-matroid}  is feasible for \eqref{eq:ar}. In particular, we have
\begin{lemma}[\bf {Feasibility}] \label{lem:feasibilty-matroid}
The affine solution in \eqref{eq:feasible-affine-matroid}  is feasible for the adjustable problem \eqref{eq:ar}.  
\end{lemma}
The proof is similar to the proof of Lemma \ref{lem:feasibilty}.

\subsection{Cost analysis.}
In the following two lemmas, we analyze the cost of the linear and static parts in our affine solution \eqref{eq:feasible-affine-matroid}.
\begin{lemma}[\bf {Cost of Linear part}] \label{lem:cost-linear-matroid}
The cost of the linear part $\mb y_{\sf Lin}(\mb h)$ defined in \eqref{eq:feasible-affine-matroid} is $ O(\beta \log n)  \cdot \opt $ for any $\mb h \in {\cal U}$.
\end{lemma}
\proof{Proof.}
We have, 
\red{
\begin{align*}
\mb d ^T\mb y_{\sf Lin}( \mb h )   =  \sum_{\ell=1}^L \sum_{i \in  {\cal I}_{\ell}} \alpha_i h_i \mb d^T {\mb v}_i =   \sum_{\ell=1}^L \sum_{i \in  {\cal I}_{\ell}} \alpha_i h_i z( \mb e_i)
& \leq    \sum_{\ell=1}^L  \beta \cdot ( {\nu}_{\ell}-   {\nu}_{\ell-1})         \sum_{i \in  {\cal I}_{\ell}}   w_{\ell i}h_i \\
& \leq   \beta  \sum_{\ell=1}^L  ( {\nu}_{\ell}-   {\nu}_{\ell-1})       \\
& =  \beta \cdot {\nu}_L  \\
& =  O(\beta \log n)  \cdot \opt.
\end{align*}
where the first inequality holds because $\alpha_i z(\mb e_i)  \leq \beta w_{\ell i} \cdot( {\nu}_{\ell}-   {\nu}_{\ell-1})  $ for all $i \in {\cal I}_{\ell}$ and $\ell \in[L]$, the second inequality holds because $ \sum_{i \in {\cal I}_{\ell}}  w_{\ell i} h_i \leq 1$ for any $ \mb h \in {\cal U}$ and the last equality follows from Lemma~\ref{lem:lastone}.
\hfill
\Halmos
\endproof
}
\begin{lemma}[\bf {Cost of Static part}] \label{lem:cost-static-matroid}
The cost of the static part $(\mb x_{\sf Sta}, \mb y_{\sf Sta})$ defined  in \eqref{eq:static-matroid} is $O( \beta\log n) \cdot \opt$.
\end{lemma}

\red{
\proof{Proof.}
Denote $\mb y_{\cal A}$ the solution provided by the online algorithm ${\cal A}$ that covers the sequence $(\mb a_1, \mb a_2, \ldots, \mb a_L)$. 
Consider
$$ {\cal J}_1 = \left\{ i \in [m] \; \bigg\vert \; ( \tilde{\mb B} \mb y_{\cal A})_i \geq \frac{1}{2}    \; \; \text{and} \; \; \alpha_i > 0 \right\}.$$
We have for all $i \in {\cal J}_1$,  $ 2 {\mb B} \mb  y_{\cal A} \geq \alpha_i \mb e_i $, i.e.,
$$ 2 {\mb B} \mb  y_{\cal A} \geq \sum_{i \in {\cal J}_1} \alpha_i \mb e_i .$$
and from Lemma \ref{lem:lastone},  
$$ 2 \mb d^T \mb y_{\cal A} = 2 {\nu}_L = O( \log n) \cdot \opt.$$ 
Now, we focus on the set of the remaining components. Denote
$${\cal J}_2= {\cal I}^c \setminus  \{  {\cal T} \cup {\cal J}_1    \}$$
and for $\ell =1, \ldots,L$ denote
$${{\cal V}}_{\ell}= {\cal J}_2 \cap S_{\ell}.$$
For each $\ell \in [L]$, we apply  the structural result in Lemma \ref{lem:structural-resulto} for the subset ${ \cal V}_{\ell}$ with parameters $ w_{\ell i} $, $ \gamma = 2 ( {\nu}_{\ell}- {\nu}_{\ell -1})$, cost vector $\mb d$ and constraint matrix $\tilde{\mb B} $. The first condition of Lemma \ref{lem:structural-resulto} is satisfied because for any $i \in { \cal V}_{\ell}$,
$$ \frac{\alpha_i z( \mb e_i )}{w_{\ell i}}  = \frac{\Theta( \mb e_i )}{w_{\ell i}}   > \beta ( {\nu}_{\ell}- {\nu}_{\ell -1}) = 4 \gamma \frac{\log n }{\log \log n}  .$$
Consider any ${\cal W} \subseteq {\cal V}_{\ell} $ such that $ \sum_{i \in {\cal W}} w_{\ell i}   \leq 1$. We have  $\mathbbm{1}({\cal W}) \in {\cal U}_{\ell}$. Moreover, $\mb y_{\cal A}$ covers less than $\frac{1}{2}\mathbbm{1}({\cal W}).$ Therefore,
$$
\frac{1}{2} \Theta\left({\cal W}\right)   \leq    {{\sf Aug} }\left( \mathbbm{1}({\cal W})  \big\vert \; \mb a_1, \mb a_2, \ldots , \mb a_L \right). \\
$$
Furthermore, 
\begin{align*}
{{\sf Aug} } \left( \mathbbm{1}({\cal W})  \big\vert \; \mb a_1, \mb a_2, \ldots , \mb a_L \right) &\leq   {{\sf Aug} }\left( \mathbbm{1}({\cal W}) \big\vert \; \mb a_1, \mb a_2, \ldots , \mb a_{\ell-1} \right) \\
&\leq { {\cal A} }\left( \mathbbm{1}({\cal W})  \big\vert \; \mb a_1, \mb a_2, \ldots , \mb a_{\ell-1} \right) \\
&\leq { {\cal A} }\left( \mb  a_{\ell}  \big\vert \; \mb a_1, \mb a_2, \ldots , \mb a_{\ell-1} \right) \\
&={\nu}_{\ell}- {\nu}_{\ell-1} , 
\end{align*}
where the first inequality holds because the  cost of covering $\mathbbm{1}({\cal W})$ given the online solution for $( \mb a_1, \mb a_2, \ldots , \mb a_{L})$ is 
smaller than the cost of covering $\mathbbm{1}({\cal W})$ given the online solution for $(\mb a_1, \mb a_2, \ldots , \mb a_{\ell-1})$. The second inequality holds because the cost of the online algorithm is less than the offline cost and the third one follows from Step 3 in the construction of the greedy scenario $ \mb a$ in Algorithm \ref{algo1}. Hence,
$$  \Theta \left( {\cal W}\right) \leq 2 ({\nu}_{\ell}- {\nu}_{\ell-1} )  = \gamma .$$
Therefore,  the second condition of Lemma \ref{lem:structural-resulto} is also satisfied and consequently,
$$ \Theta({\cal V}_{\ell}) \leq  4 \gamma \frac{\log n}{\log \log n}  =        \beta   ( {\nu}_{\ell}- {\nu}_{\ell-1}).$$
By taking the sum over all $\ell=1,\ldots,L$, we get
$$ \sum_{\ell=1}^L \Theta({\cal V}_{\ell}) \leq \beta {\nu}_L = O(\beta  \log n ) \cdot  \opt.$$
For $\ell \in [L]$, denote $\mb y_{\ell}$ an optimal solution corresponding to $\Theta({\cal V}_{\ell})$. In particular,
$$ \sum_{\ell=1}^L \mb B \mb y_{\ell} \geq \sum_{i \in {\cal J}_2} \alpha_i \mb e_i.$$
By feasibility of the optimal solution, we have
$$ \mb A \mb x^*+ \mb B \mb y^*(\mb 0) \geq \mb 0,$$
i.e.,
$$ \mb B \mb y^*(\mb 0) \geq \sum_{ i=1}^m (\alpha_i -1) \mb e_i.$$
Moreover, since $\mb B$ and $y^*(\mb 0) $ are non-negative, we have
$$ \mb B \mb y^*(\mb 0) \geq \sum_{ i=1}^m (\alpha_i -1)^+ \mb e_i.$$
Therefore, we have the following candidate static solution for~\eqref{eq:static-matroid}:
\begin{align*}
\mb x_{\sf Sta} & = \mb x^* \\
\mb y_{\sf Sta} & = \mb y^*(\mb 0) + 2 \mb y_{\cal A} + \sum_{\ell =1}^L \mb y_{\ell}.
\end{align*}
Putting it all together,
\begin{align*}
\mb A \mb x^* + \mb B\mb y^*(\mb 0)+ 2\mb B \mb y_{\cal A}+ \sum_{\ell=1}^L \mb B \mb y_{\ell} &\geq \sum_{i=1}^m (1-\alpha_i) \mb e_i + \sum_{ i=1}^m (\alpha_i -1)^+ \mb e_i + \sum_{i \in {\cal J}_{1} } \alpha_i \mb e_i+\sum_{i \in {\cal J}_{2} } \alpha_i \mb e_i \\
&= \sum_{i=1}^m (1-\alpha_i)^+ \mb e_i  + \sum_{i \in {\cal J}_{1} \cup  {\cal J}_{2}  } \alpha_i \mb e_i \\
&= \sum_{i \in {\cal I}}   (1-\alpha_i)^+ \mb e_i +\sum_{i \in {\cal I}^c}   (1-\alpha_i)^+ \mb e_i + \sum_{i \in {\cal J}_{1} \cup  {\cal J}_{2}  } \alpha_i \mb e_i \\
&=    \sum_{i \in {\cal I}}     (1-\alpha_i)^+ \mb e_i+ \sum_{i \in {\cal \tau}}    (1-\alpha_i)^+ \mb e_i  +  \sum_{i \in {\cal J_1}  \cup {\cal J_2} }     (1-\alpha_i)^+ \mb e_i + \sum_{i \in {\cal J}_{1} \cup  {\cal J}_{2}  } \alpha_i \mb e_i \\
& \geq    \sum_{i \in {\cal I}}     (1-\alpha_i)^+ \mb e_i+ \sum_{i \in {\cal \tau}}     \mb e_i  +  \sum_{i \in {\cal J_1}  \cup {\cal J_2} }      \mb e_i \\
& =  \sum_{i \in {\cal I}}    (1-\alpha_i)^+ \mb e_i  +   \sum_{i \in {\cal I}^c}      \mb e_i
\end{align*}
where the last inequality holds because $ \alpha_i \leq 0 $ for all $ i \in {\cal \tau}$ and $   (1-\alpha_i)^+ + \alpha_i \geq 1$ for all $ i \in  {\cal J}_1 \cup {\cal J}_2$.
Therefore,
\begin{align*}
z_{\sf Sta} &\leq  \mb c^T \mb x^*+ \mb d^T \mb y^*(\mb 0)+ 2 \mb d^T \mb y_{\cal A} + \sum_{\ell=1}^L z({\cal V}_{\ell})   \\
&\leq  \opt+ O(\log n ) \cdot \opt+ O( \beta \log n )  \cdot \opt \\
& = O( \beta\log n) \cdot \opt.
\end{align*} 
\hfill
\Halmos
\endproof
}

\proof{Proof of Theorem \ref{thm:partition-matroid}.}
Lemma \ref{lem:feasibilty-matroid} show that our affine solution \eqref{eq:feasible-affine-matroid} is feasible for the adjustable problem \eqref{eq:ar}. Lemma \ref{lem:cost-linear-matroid} and \ref{lem:cost-static-matroid} show that  the cost of the feasible affine solution is less than
$$O(\beta \log n)\cdot \opt +O(\beta \log n)\cdot \opt =O \left( \frac{\log^2 n}{\log \log n} \right) \cdot \opt $$ which implies that, $$ z_{\sf Aff}({\cal U}) = O \left( \frac{\log^2 n}{\log \log n} \right)\cdot z_{\sf AR}({\cal U}).$$
\hfill
\Halmos
\endproof

We would like to note that Gupta et al. \cite{gupta2016robust} give $O(\log n)$-approximation to \eqref{eq:ar} in the special case $\mb A, \mb B \in \{0,1\}^{m \times n}$, $ \mb d = \lambda \mb c$ and $w_{\ell i} = w$ are all $\ell, i$ for some constant $w$. Therefore, for this special case the bound of~\cite{gupta2016robust} is stronger than our bound in Theorem \ref{thm:partition-matroid}. However, their algorithm does not give a functional policy approximation. Here, our focus is different, namely, to analyze the performance of affine policies that are widely used in practice and exhibit strong empirical performance. Our analysis shows that the performance of affine policies for disjoint constrained budgeted sets is near-optimal and nearly matches the hardness of the problem. Note that our bound in Theorem~\ref{thm:partition-matroid} is not necessarily tight. It is an interesting open question to study if affine policies also give an optimal approximation for this more general class of budgeted uncertainty sets.

\section{General Intersection of budgeted  uncertainty sets.}
In this section, we consider the general intersection of budget of uncertainty sets given by \eqref{eq:intersection:set}
$${\cal U} = \left\{ \mb h \in [0,1]^m  \; \Big\vert \;   \sum_{ i \in S_{\ell}} w_{\ell i} h_i \leq 1 \;\;  \forall \ell \in [L] \right\},$$
where $  \mb w_{\ell} \in [0,1]^m$ and $S_{\ell}$ for $\ell \in [L]$ is a general family of subsets of $[m]$. This class is a generalization of the single budget of uncertainty set \eqref{eq:budget}. It captures many important sets including CLT sets considered in Bertsimas and Bandi \cite{bandi2012tractable} and {\em inclusion-constrained budgeted} sets considered in Gounaris et al. \cite{gounaris2014adaptive}. In this section , we study the performance of affine policies for intersection of budget of uncertainty sets \eqref{eq:intersection:set} and show strong theoretical guarantees. We start by the case when the set \eqref{eq:intersection:set} verifies some symmetric properties ({\em permutation invariant sets}) and then we give our results for the general form \eqref{eq:intersection:set}.

\subsection{Permutation Invariant Sets.}
We consider intersection of budgeted  sets that are \textit{permutation invariant}. 
\begin{definition}[\bf {Permutation Invariant Sets}] \label{def:perm}
We say that ${\cal U}$ is a  \textit{permutation invariant set} if $\mb{x}\in{\cal U}$ implies that for any permutation $\tau$ of $\{1,2,\ldots,m\}$, $\mb{x}^{\tau}\in{\cal U}$ where $x^\tau_i= x_{\tau(i)}$.
\end{definition}
This class of sets captures many important sets including CLT sets that have been considred in Bertsimas and Bandi \cite{bandi2012tractable}. Note that a CLT set is given by
\begin{equation} \label{CLT:sets}
{\cal U} = \left\{ \mb h \in [0,1]^m  \; \Big\vert \;   \sum_{i \in {\cal S}} h_i \leq \Gamma  \;\;  \forall { \cal S} \subseteq [m], \vert {\cal S} \vert =k \right\},
\end{equation}
for some $k \in \mathbb{N}$. The following theorem gives our performance bound for affine policies under the class of intersection of budgeted sets that are permutation invariant.

\begin{theorem}\label{thm:permutation-inv}
Consider the two-stage adjustable problem \eqref{eq:ar} where $ {\cal U}$ is the intersection of $L$ budget constraints \eqref{eq:intersection:set}.  Suppose that $\cal U$ is permutation invariant set \red{and $\cal X$ is a polyhedral cone.}  Then,
$$ z_{\sf Aff}({\cal U}) = O \left( \log L \cdot \frac{\log n}{\log \log n} \right) \cdot z_{\sf AR}({\cal U}) .$$
\end{theorem}

\proof{Proof.}
Our proof relies on a geometric property that we show for budgeted uncertainty sets that are permutation invariant. In particular, we show that  for any $\cal U$ permutation invariant, there exists a (single) budget of uncertainty set $\cal V$ of the form \eqref{eq:budget} such that
\begin{equation}\label{inclusion}
 \frac{1}{4 \log L} \cdot {\cal V} \subseteq {\cal U}  \subseteq  2  {\cal V}.
\end{equation}
Since ${\cal U}$ is permutation invariant, 
$$ \gamma  \mb e \in \red{{\argmax}} \left\{ \mb e^T \mb h \; \bigg\vert \; \mb h \in {\cal U} \right\},$$
for some $ \gamma \in [0,1]$. Consider
$$  {\xi}_i \overset{i.i.d.}{\sim} {\sf Ber} ( \gamma) \qquad  i=1,\ldots,m,$$
i.e., ${\xi}_1, {\xi}_2, \ldots, {\xi}_m$ are i.i.d. Bernoulli random variables of parameter $ \gamma$. Let $$\mb {\xi} = ({\xi}_1, {\xi}_2, \ldots, {\xi}_m).$$ 
Consider the following  budget of uncertainty set 
$$ \tilde{\cal V} = \left\{ \mb h \in [0,1]^m  \; \Bigg\vert \;   \sum_{i=1}^m h_i \leq  \sum_{i=1}^m {\xi}_i  \right\}.$$
Note that  $\tilde{\cal V}$ is random depending on the realization of ${\xi}_1, {\xi}_2, \ldots, {\xi}_m$ . We show that 
$$\mathbb{P} \left( \frac{1}{4 \log L} \cdot \tilde{\cal V} \subseteq {\cal U}  \subseteq  2  \tilde{\cal V} \right) > \epsilon ,$$ for some constant $\epsilon >0$ which implies the existence of $\tilde{\cal V}$  such that \eqref{inclusion} is verified. For that purpose, we show first that the right inclusion holds with a constant probability and then the left one holds with high probability.

\begin{claim}\label{claim:5.1} 
$\mathbb{P} \left( {\cal U} \subseteq 2  \tilde{\cal V} \right) \geq  1-  e^{-\frac{1}{8}}.$
\end{claim}

\noindent
Let us prove the above claim. Note that,
 $$ \gamma  m =\underset{\mb h \in {\cal U}}{\max} \;  \mb e^T \mb h. $$
 Suppose that, $$\gamma  m \leq 2 \mb e^T \mb {\xi}.$$ Then for all $\mb h \in{\cal U}$,
 $$\mb e^T \mb h \leq 2 \mb e^T \mb {\xi}$$
i.e., for all $\mb h \in{\cal U}$
 $$ \mb h \in 2 \tilde{\cal V}.$$ 
Hence, $\gamma  m \leq 2 \mb e^T \mb {\xi}$ implies that ${\cal U} \subseteq 2  \tilde{\cal V}$.
Therefore,
\begin{align*}
 \mathbb{P} \left( {\cal U} \subseteq 2  \tilde{\cal V} \right)  & \geq  \mathbb{P} \left(    2 \mb e^T \mb{\xi} \geq \gamma  m   \right) \\
 &=   \mathbb{P} \left( \sum_{i=1}^m {\xi}_i \geq \frac{1}{2} \gamma m   \right).
\end{align*}
We know that $ \mathbb{E} \left( \sum_{i=1}^m {\xi}_i \right) = \gamma m $. Therefore, from the Chernoff inequality in Lemma \ref{lem:chernof2}, we have,
$$ \mathbb{P} \left( \sum_{i=1}^m {\xi}_i \geq \frac{1}{2} \gamma m   \right) \geq  1- \exp \left(- \frac{\gamma m}{8} \right) \geq   1-  e^{-\frac{1}{8}}.$$
where the last inequality holds because $\gamma m \geq 1$ since $\mb e_i \in {\cal U}$ for all $i$ and $\cal U$ is convex.

\begin{claim}\label{claim:5.2} 
$\mathbb{P} \left(  \tilde{\cal V} \subseteq  4 \log L \cdot {\cal U} \right) \geq 1- \frac{1}{L}.$
\end{claim}

\noindent
Note that $\mb {\xi}$ is an extreme point of $\tilde{\cal V}$ and that all pareto extreme points of $\tilde{\cal V}$  are just permutation of $\mb  {\xi}$. Moreover, we know that $\cal U$ is permutation invariant set, hence  if $\cal U$ contains $ \mb {\xi}$ then $\cal U$ contains all pareto extreme points of $\cal U$ and consequently contains $\cal U$ by down-monotonicity. Therefore,
\begin{align*}
\mathbb{P} \left(  \tilde{\cal V} \subseteq  4 \log L \cdot {\cal U} \right)  & \geq  \mathbb{P} \left(  \mb {\xi} \in  4 \log L \cdot {\cal U} \right)  \\
& =  \mathbb{P} \left(  \mb {\nu}_{\ell}^T \mb {\xi} \leq   4 \log L, \; \; \forall \ell \in  [L] \right)  \\
& =  1-  \mathbb{P} \left(  \exists \ell \in [L],      \;   \mb {\nu}_{\ell}^T \mb {\xi} >   4 \log L   \right)  \\
& \geq  1-  \sum_{\ell=1}^L \mathbb{P} \left(   \mb {\nu}_{\ell}^T \mb {\xi} >   4 \log L   \right),
\end{align*}
where the last inequality follows from a union bound. We have $ \mathbb{E}( \mb {\nu}_{\ell}^T \mb {\xi} ) = \mb {\nu}_{\ell}^T \gamma \mb e \leq 1 \leq \log L$ \red{for $L \geq 2$} because $\gamma \mb e $ is a feasible point in $\cal U$.  Therefore, from Lemma \ref{lem:chernof1} with $\delta=3$.
$$ \mathbb{P} \left(   \mb {\nu}_{\ell}^T \mb {\xi} >   4 \log L   \right) \leq \left( \frac{e^3}{4^4}\right)^{\log L} \leq \left(e^{-2} \right)^{\log L}  = \frac{1}{L^2}.$$
We conclude that,
\begin{align*}
\mathbb{P} \left(  \tilde{\cal V} \subseteq  4 \log L \cdot {\cal U} \right)  \geq  1-  \sum_{\ell=1}^L  \frac{1}{L^2} = 1- \frac{1}{L} . \\
\end{align*}
Hence, from Claim~\ref{claim:5.1} and Claim~\ref{claim:5.2}, there exists a budget of uncertainty set $\tilde{\cal V}$ with a non zero probability  that verifies the inclusion in \eqref{inclusion}.
Therefore,
\begin{align*}
z_{\sf Aff}({\cal U}) & \leq 2 \cdot   z_{\sf Aff}( \tilde{\cal V}) \\
&= 2 \cdot O \left( \frac{\log n}{\log \log n} \right) \cdot   z_{\sf AR}( \tilde{\cal V}) \\
&\leq 2 \cdot O \left( \frac{\log n}{\log \log n} \right)  \cdot 4 \log L \cdot   z_{\sf AR}( {\cal U})= O\left( \frac{\log n}{\log \log n}\ \cdot \log L \right) \cdot z_{\sf AR}({\cal U}),
\end{align*}
where the first inequality holds because $ {\cal U} \subseteq  2 \tilde{\cal V}$ \red{and $ 2 \cdot\cal X \subseteq  {\cal X}$ ($\cal X$ is a polyhedral cone).} The first equality follows from Theorem \ref{thm:budget} because $\tilde{\cal V}$ is a budget of uncertainty set, and finally the last inequality holds because $  \tilde{\cal V} \subseteq  4 \log L \cdot     \cal U$ \red{and  $ 4 \log L \cdot\cal X \subseteq  {\cal X}$ ($\cal X$ is a polyhedral cone).}
\hfill
\Halmos
\endproof

We would like to note that the result of Theorem \ref{thm:permutation-inv} extends as well to the class of  intersection of budgeted  sets that are {\em scaled permutation invariant}. We say that ${\cal U} $ is a scaled permutation invariant set if there exists $\mb \lambda \in \mathbb{R}^m_+$ and $ {\cal V}$ a permutation invariant set such that 
$$ {\cal U}= {\sf diag} (\mb \lambda) \cdot {\cal V}.$$ 
In fact, for a given scaled permutation invariant set ${\cal U}$, it is possible to scale the two-stage adjustable problem \eqref{eq:ar} and get a new problem where the uncertainty set is permutation invariant. Indeed, suppose $ {\cal U}= {\sf diag} (\mb \lambda) \cdot {\cal V}$ where $ {\cal V}$ is a permutation invariant set; by multiplying the constraint matrices $\mb A$ and $ \mb B$ by ${\sf diag} (\mb \lambda)^{-1}$, we get a new problem where the uncertainty set now is permutation invariant. The performance of affine policy is not affected by this scaling and the bound given by Theorem \ref{thm:permutation-inv} still hold.

\subsection{General Intersection of Budgets.}
Consider the general intersection of budgeted  sets  \eqref{eq:intersection:set} which is given by
$${\cal U} = \left\{ \mb h \in [0,1]^m  \; \Big\vert \; \sum_{i \in S_{\ell}}  w_{\ell i}  h_i \leq 1 \;\;  \forall \ell \in [L] \right\}.$$
We show that affine policy gives a worst-case bound of $O \left( L \cdot  \frac{\log n}{\log \log n} \right) $ where $L$ is the number of constraints in  $\cal U$. In particular, we have the following theorem.
\begin{theorem}\label{thm:general-intersection}
Consider the two-stage adjustable problem \eqref{eq:ar} where $ {\cal U}$ is the  intersection of $L$ budgeted  sets given by \eqref{eq:intersection:set}. \red{ Suppose that $\cal X$ is a polyhedral cone.} Then,
$$ z_{\sf Aff}({\cal U}) = O \left( L \cdot  \frac{\log n}{\log \log n} \right)    \cdot z_{\sf AR}({\cal U}) .$$
\end{theorem}

\proof{Proof.}
Denote for all $\ell \in [L]$,
$$\mb w_{\ell} = \sum_{i \in S_{\ell}} w_{\ell i} \mb e_i$$
and let
$$  \mb{\bar w} =  \frac{1}{L}\sum_{\ell=1}^L \mb w_{\ell} .$$
Consider the following budget of uncertainty set,
$${\cal V} = \left\{ \mb h \in [0,1]^m  \; \Bigg\vert \;  \mb{\bar w}^T \mb h \leq 1 \right\}.$$

We show that ${\cal U }  \subseteq {\cal V} \subseteq  L \cdot {\cal U }$.  Suppose $ \mb h \in {\cal U}$. Then, for any ${\ell} \in [L]$, we have $ \mb w_{\ell}^T \mb h \leq 1$. Therefore, by summing up all these inequalities and dividing by $L$, we get $\mb{\bar w}^T \mb h \leq 1$, i.e., $\mb h \in {\cal V}$. Hence ${\cal U }  \subseteq {\cal V}$. Conversely, suppose $ \mb h \in {\cal V}$. For any ${\ell} \in [L]$, we have $ \mb w_{\ell}^T \mb h \leq  \sum_{\ell=1}^L   \mb w_{\ell}^T \mb h   \leq L$, hence $\mb h \in L \cdot {\cal U}$ and consequently ${\cal V} \subseteq  L \cdot {\cal U }$. Therefore,
\begin{align*}
z_{\sf Aff}({\cal U}) & \leq    z_{\sf Aff}( {\cal V}) \\
&= O \left( \frac{\log n}{\log \log n} \right) \cdot   z_{\sf AR}( {\cal V}) \\
& \leq  O\left( \frac{\log n}{\log \log n} \right) \cdot L \cdot z_{\sf AR}({\cal U}),
\end{align*}
where the first inequality holds because $ {\cal U} \subseteq {\cal V}$, the second one is a consequence of Theorem \ref{thm:budget} since ${\cal V}$ is a budget of uncertainty set of the form \eqref{eq:budget}, and finally the last inequality holds because $ {\cal V} \subseteq   L \cdot     {\cal U}$ \red{and  $ L \cdot\cal X \subseteq  {\cal X}$ ($\cal X$ is a polyhedral cone).}

\hfill
\Halmos
\endproof


\section{\red{Faster algorithm for near-optimal affine solutions.}} \label{sec:numerical}

In this section, we present an algorithm to compute an approximate affine policy for~\eqref{eq:ar} under budget of uncertainty sets, that is significantly faster than solving the optimization program  \eqref{eq:aff:LP} that computes the optimal affine policy. Our algorithm is based on the analysis of the performance of affine policies  that shows the existence of a good affine solution that satisfies certain nice structural properties. In particular, our construction of approximate affine solution in Section 3.1 partitions the components into {\em expensive} and {\em inexpensive} components based on a threshold. We cover \red{ a fraction of the {\em inexpensive} components using a linear solution and the remaining components using a static solution.}  In particular, we show that there exists an affine solution with such a structure and cost at most $O \left( \log n /\log \log n \right)$ times the optimal optimal cost of \eqref{eq:ar} for some partition of components into {\em expensive} and {\em inexpensive}. Based on this structure, we give a faster algorithm to compute an approximate affine solution for budget of uncertainty sets. 

\red{
\subsection{Our algorithm.}
Let $\cal U$ be the budget of uncertainty set \eqref{eq:budget} given by
$$ {\cal U} = \left\{ \mb h \in [0,1]^m  \; \Big\vert \;   \sum_{i=1}^m w_i h_i \leq 1 \right\}.$$
Recall from the proof of Theorem  \ref{thm:budget}, we construct our candidate affine solution by partitioning the components of $[m]$ into two subsets $ {\cal  I}$ and its complement $ {\cal  I}^c$.  The linear solution \eqref{eq:linear-sol} is given by 
\begin{equation*} \label{eq:linear-sol-algo}
 \mb y_{\sf Lin}( \mb h )   = \sum_{i \in {\cal I} } \alpha_i h_i {\mb v}_i 
\end{equation*}
 where
$${\cal  I}= \left\{ i \in [m] \; \bigg\vert \;  \frac{\alpha_i z( \mb e_i)}{w_i}   \leq \beta \cdot \opt \right\},$$
and for all $i \in [m]$,
\[
\begin{aligned}
\alpha_i & = 1 -( \mb A \mb x^*)_i, \\
\mb v_i & \in  \argmin_{\mb y \geq \mb 0} \left\{ \mb d^T \mb y \; \bigg\vert \;   \mb B \mb y \geq  \mb e_i  \right\}.
\end{aligned}
\]
Let 
$$ \mb Y = [ \mb v_1 \; \vert \; \mb v_2 \; \vert \ldots \vert \mb v_m ].$$
Based on the structure of the linear part, we propose the following approximate affine solution: 
$$ \mb y ( \mb h )   =  \mb Y \cdot {\sf diag}  ( \mb {\alpha} ) \cdot \mb h + \mb q $$
where  $\mb Y$ is a constant that can be computed efficiently upfront and $   {\alpha}_i $ for $i \in [m]$ are non-negative variables. This structure captures our candidate solution  \eqref{eq:linear-sol}. 
Hence, we reduce the number of second stage variables  from $O(nm)$ in~\eqref{eq:aff} to $O(n+m)$. Moreover, the non-negativity constraint on $ \mb y ( \mb h )$ reduces to $\mb \alpha \geq \mb 0$ and $\mb q \geq \mb 0$ in this special class of affine solutions. Restricting to the above class of affine solutions, we have the following optimization problem.
\begin{equation}\label{eq:aff-approx}
\begin{aligned}
 \min_{\mb x, \mb \alpha, \mb q} \; & \mb{c}^T \mb{x} + \max_{\mb{h}\in {\cal U}} \; \mb{d}^T\left( \mb Y \cdot {\sf diag}  ( \mb {\alpha} ) \cdot \mb{h}+\mb{q}\right) \\
& \mb{A}\mb{x} + \mb{B}\left( \mb Y \cdot {\sf diag}  ( \mb {\alpha} )\cdot \mb{h}+\mb{q}\right) \; \geq \; \mb{h},   \; \; \; \forall \mb h \in {\cal U} \\
& \red{\mb{x}  \in  {\cal X} }, \mb{\alpha}   \in  {\mathbb R}^{m}_+, \mb{q}   \in  {\mathbb R}^{n}_+.\\
\end{aligned}
\end{equation}
Using similar reformulations as in Lemma \eqref{lem:LP-form}, the above problem can be formulated as the following LP:
\begin{equation}\label{eq:aff:approx2}
\begin{aligned}
 \min \; &  \mb{c}^T \mb{x} + z \\
& z- \mb d^T \mb q \geq \mb r^T \mb v \\
& \mb R^T \mb v \geq \mb Y \cdot {\sf diag}  ( \mb {\alpha} )^T \mb d \\
& \mb A \mb x + \mb B \mb q \geq \mb V^T \mb r\\
& \mb R^T \mb V \geq \mb I_m - \mb B \mb Y \cdot {\sf diag}  ( \mb {\alpha} ) \\
& \red{\mb{x}  \in  {\cal X} }, \; \mb{v}   \in  {\mathbb R}^{L}_+, \; \mb{U}   \in  {\mathbb R}^{L \times n}_+, \; \mb{V}   \in  {\mathbb R}^{L \times m}_+ \\
&  \mb{\alpha}   \in  {\mathbb R}^{m}_+ \; \mb q \in \mathbb{R}_+^n, \; z \in \mathbb{R}.   \\
\end{aligned}
\end{equation}
The above formulation is significantly faster than solving~\eqref{eq:aff} as we observe in our numerical experiments. Algorithm~\ref{algo2} describes the detail of our algorithm.
\begin{algorithm}[H]
\caption{Computing Approximate Affine Policy }\label{algo2}
\begin{algorithmic}[1]
 \For{ $i=1,\ldots,m$}
 \State $$ \mb v_i \in  \argmin_{\mb y \geq \mb 0} \left\{ \mb d^T \mb y \; \bigg\vert \;   \mb B \mb y \geq  \mb e_i  \right\}.$$
 \EndFor
 \State  $$ \mb Y = [ \mb v_1 \; \vert \; \mb v_2 \; \vert \ldots \vert \mb v_m ].$$
 \State Solve the LP :
 \begin{equation*}\label{eq:aff:approx22}
\begin{aligned}
z_{\sf Alg}= \min \; &  \mb{c}^T \mb{x} + z \\
& z- \mb d^T \mb q \geq \mb r^T \mb v \\
& \mb R^T \mb v \geq \mb Y \cdot {\sf diag}  ( \mb {\alpha} )^T \mb d \\
& \mb A \mb x + \mb B \mb q \geq \mb V^T \mb r\\
& \mb R^T \mb V \geq \mb I_m - \mb B \mb Y \cdot {\sf diag}  ( \mb {\alpha} ) \\
& \red{\mb{x}  \in  {\cal X} }, \; \mb{v}   \in  {\mathbb R}^{L}_+, \; \mb{U}   \in  {\mathbb R}^{L \times n}_+, \; \mb{V}   \in  {\mathbb R}^{L \times m}_+ \\
&  \mb{\alpha}   \in  {\mathbb R}^{m}_+ \; \mb q \in \mathbb{R}_+^n, \; z \in \mathbb{R}.   \\
\end{aligned}
\end{equation*}
\State \Return $z_{\sf Alg}.$
\end{algorithmic}
\end{algorithm}
We would like to note that since our approximate affine solution is based on the construction of affine policy in our analysis, the worst-case approximation bound for our approximate affine solution  is also $O (  \frac{\log n}{\log \log n} )$.
}

\subsection{Numerical experiments.}

We study the empirical performance of our algorithm for budget of uncertainty sets both from the perspective of computation time and the quality of the solution. 

\vspace{2mm}
\noindent
{\bf Experimental setup.}
We use the same test instances as in Ben-Tal et al. \cite{ben2018tractable}. In particular, we choose $n=m$, $\mb c= \mb d = \mb e $ and $ \mb A= \mb B$ where $\mb B$ is  randomly generated as $ \mb B = \mb I_m + \mb G,$ where $\mb I_m$ is the identity matrix and $\mb G$ is a random normalized Gaussian, i.e. $G_{ij} = \vert Y_{ij} \vert / \sqrt{m}$ where $Y_{ij}$ are i.i.d. standard Gaussians.
 Let consider  the following budget of uncertainty sets:
\begin{eqnarray}
\label{set:k-ones:num} {\cal U}_1 & = & \left\{ \mb h \in [0,1]^{m} \;  \Bigg\vert \; \sum_{i=1}^m h_i \leq k \right\} \\
\nonumber & \\
\nonumber & \\
\label{set:budget:num}  {\cal U}_2 & = & \left\{ \mb h \in [0,1]^{m} \;  \Bigg\vert \; \sum_{i=1}^m w_i h_i \leq 1 \right\}.
\end{eqnarray}

For our numerical experiments, we choose $k=c  \sqrt{m}$  with $c$  a random uniform constant between $1$ and $2$ for the first uncertainty set ${\cal U}_1$. For the second uncertainty set ${\cal U}_2$, we choose $\mb w$ a normalized Gaussian vector, i.e., $ w_i = \vert G_i \vert / \|  {\mb G} \|_2 $ where $G_i$ are i.i.d. standard Gaussians. We consider values of $m$ from $m=10$ to $ m=200$ in increments of $10$ and consider $20$ instances for each value of $m$. 

We compute the optimal affine solution by solving the LP formulation \eqref{eq:aff:LP}.  We compute our approximate affine solution returned by Algorithm \ref{algo2}. We denote  $z_{\sf{Alg}} ({\cal U})$ and $ z_{\sf{Aff}}({\cal U})$ respectively the cost of our affine solution returned by Algorithm \ref{algo2} and the cost of the optimal affine solution. For each $m$ from $m=10$ to $m=200$, we report the average ratio $  z_{\sf{Alg}} ({\cal U})/ z_{\sf{Aff}}({\cal U}) $, the running time of Algorithm \ref{algo2} in seconds $(T_{\sf Alg}(s))$ and the running time of the optimal affine policy in seconds $(T_{\sf aff}(s))$. We present the results of our computational experiments in Table \ref{tab:num}. The numerical results are obtained using Gurobi 7.0.2 on a 16-core server with 2.93GHz processor and 56GB RAM.

\vspace{2mm}
\noindent
{\bf Results.}
\red{ We observe from Table \ref{tab:num} that our algorithm is significantly faster than the optimal affine policy. In fact, Algorithm \ref{algo2} scales very well and the average running time is only a few seconds even for large values of $m$. On the other hand, computing the optimal affine solution becomes computationally challenging for large values of $m$.  For example,  for $m=100$, the average running time is around 3 minutes for ${\cal U}_1$ and more than 11 minutes for ${\cal U}_2$. For $m=200$, the average running time is more than an hour for ${\cal U}_1$ and more than 3 hours for ${\cal U}_2$.  Furthermore, we observe that the gap between our affine solution and the optimal one is under $15\%$. Moreover, this gap does not increase with the dimension of $m$, thereby confirming that our affine solution performs well even for large values of $m$.}

\begin{table}[htp]
\begin{subtable}{0.5\linewidth}\centering
{\begin{tabular}{|c|c|c|c|}
\hline
$m$ & $T_{\sf aff}(s)$ & $T_{\sf Alg}(s)$ & $z_{\sf{Alg}}/ z_{\sf{Aff}} $ \\ \hline
10  & 0.009	&	0.022	&	1.146   \\\hline
20  & 0.137	&	0.105	&	1.111   \\\hline
30  &0.304	&	0.300	&	1.155   \\\hline
40  &1.268	&	0.692	&	1.126   \\\hline
50  &4.007	&	1.370	&	1.120   \\\hline
60  &9.461	&	3.089	&	1.135    \\\hline
70  &17.38	&	3.417	&	1.147   \\\hline
80  &44.75  &	5.626	&	1.103    \\\hline
90  &80.20	&	10.18	&	1.114   \\\hline
100  &153.3	&	13.23	&	1.149    \\\hline
200  &5137	&	69.33	&	1.061    \\\hline
\end{tabular}}
\caption{Budget of uncertainty  \eqref{set:k-ones:num}}\label{tab:U1}
\end{subtable}%
\begin{subtable}{.5\linewidth}\centering
{\begin{tabular}{|c|c|c|c|}
\hline
$m$ & $T_{\sf aff}(s)$ & $T_{\sf Alg}(s)$ &$z_{\sf{Alg}}/ z_{\sf{Aff}} $ \\ \hline
10  & 0.011 &		0.021	&	1.108   \\\hline
20  &0.200	&	0.110	&	1.092     \\\hline
30  &1.219	&0.353  &	1.103    \\\hline
40  &4.887	&	0.812	&	1.093   \\\hline
50  &17.13	&	1.388	&	1.096   \\\hline
60  &54.03	&	2.259	&	1.086   \\\hline
70  &129.7  	&3.625	&	1.088   \\\hline
80  &248.1 &		5.069	&	1.082  \\\hline
90  &390.9 &		6.381	&	1.080  \\\hline
100  &692.9& 	8.705	&	1.082 \\\hline
200  &**& 	68.62	&	** \\\hline
\end{tabular}}
\caption{Budget of uncertainty \eqref{set:budget:num}}\label{tab:U2}
\end{subtable}
\caption{Comparison on the performance and computation time of the optimal affine policy and our approximate affine policy. For 20 instances, we compute $ z_{\sf{Alg}} ({\cal U})  / z_{\sf{Aff}}  ({\cal U})       $ for $\cal U$ the  budget of uncertainty sets \eqref{set:k-ones:num} and \eqref{set:budget:num}. Here, $T_{\sf Alg}(s)$ denotes the running time for our approximate affine policy and $T_{\sf aff}(s)$ denotes the running time for affine policy in seconds. $**$ denotes the cases when we set a time limit of 3 hours. These results are obtained using Gurobi 7.0.2 on a 16-core server with 2.93GHz processor and 56GB RAM.}
\label{tab:num}
\end{table}

\red{
\section{General case of recourse matrix $\mb B$.} \label{section:general B}
In this section, we consider the two-stage adjustable problem $\eqref{eq:ar}$ with general recourse matrix $\mb B$ where we relax the non-negativity assumption on $\mb B$. In particular, we consider cases where some of the coefficients in $B$ could be negative. 
We show that in this case, the gap between the optimal affine policy given by \eqref{eq:aff} and the optimal adjustable problem \eqref{eq:ar} could be arbitrary large. Therefore, the non-negativity assumption on the coefficients of $\mb B $ is crucial for affine policies to have a good performance with respect to the optimal adjustable solution. 
\subsection{Large Gap Instances.}
We consider a two-stage lot-sizing problem to construct a family of instances of~\eqref{eq:ar} with general recourse matrix $\mb B$ such that the gap between the optimal adjustable solution and optimal affine policies is unbounded.}

\vspace{2mm}

\red{\noindent{\bf Two-Stage Robust Lot-Sizing Problem.} We are given a set of $m$ nodes with pairwise distances $d_{ij}$ between node $i$ and node $j$. Each node $i \in [m]$ has cost $c_i$ per unit inventory at node $i$ and has a capacity of $K_i$. Each node $i$ faces an uncertain demand $h_i$ that is realized in the second-stage. In the first-stage, the decision maker needs to decide the inventory levels, $x_i$ for each node $i \in [m]$. We model uncertain demand as an adversarial selection from a pre-specified uncertainty set ${\cal U}$ after the adversary observes the first-stage inventory decisions. In the second-stage, the decision maker can make recourse transportation decisions after observing the uncertain demand to satisfy it using the first-stage inventory. The goal is to make the first-stage inventory decisions such that the sum of first-stage inventory costs and the worst case second-stage transportation costs is minimized. This problem has been studied extensively in the literature (see for example Bertsimas and de Ruiter \cite{bertsimas2016duality}).}

\red{We can formulate the above problem in our framework of~\eqref{eq:ar} where the recourse matrix $\mb B$ is a network matrix with entries in $\{-1,0,1\}$. 
The epigraph formulation is the following.
\begin{equation}\label{eq:lot-sizing} 
\begin{aligned}
\min_{\mb x, z} \; &  \sum_{i=1}^m c_i x_i + z \\
& z \geq \sum_{i=1}^m \sum_{j=1}^m  d_{ij} y_{ij} ( \mb h) \\
&  x_i+  \sum_{j=1}^m y_{ji}( \mb h) - \sum_{j=1}^m y_{ij} ( \mb h) \geq h_i , \; \; \; \forall i \in [m], \; \; \forall \mb h \in {\cal U} \\
&  0 \leq x_i \leq K_i ,\; \; \; \forall i \in [m]\\
& \mb{y}(\mb{h}) \in  {\mathbb R}^{m^2}_+, \;  \; \; \forall \mb h \in {\cal U}.
\end{aligned}
\end{equation}}


\vspace{2mm}

\red{\noindent{\bf Family of instances.}
We consider the following family of instances for the robust lot-sizing problem~\eqref{eq:lot-sizing}. Consider a bipartite network $(J_1, J_2)$ where $ \vert J_1 \vert = \vert J_2 \vert = m/2$ ($m$ is even). We consider a budget of uncertainty set to model demand uncertainty. The inventory cost $c_i$, capacity $K_i$ for all $i\in [m]$, distances $d_{ij}$, $i,j \in [m]$ and the formulation for the uncertainty set are given as follows.
\begin{equation}\label{eq:lot-sizing-instance}
\begin{aligned}
c_i&= \left\{
    \begin{array}{ll}
        0 & \qquad\mbox{if } i \in J_1 \\
        1 & \qquad \mbox{if } i \in J_2
    \end{array}
\right. \\
& \\
K_i &=1 \qquad \forall  i \in [m] \\
& \\
d_{ij}&= \left\{
    \begin{array}{ll}
        0 & \qquad\mbox{if } i \in J_1, j \in J_2 \\
        \infty & \qquad \mbox{otherwise}.
    \end{array}
\right. \\
& \\
{\cal U} & = \left\{ \mb h \in [0,1]^{m} \;  \Bigg\vert \; \sum_{i=1}^m h_i \leq m/2 \right\}.
 \end{aligned}
 \end{equation}
For the above family of instances \eqref{eq:lot-sizing-instance}, we show that the gap between optimal affine and adjustable policies is unbounded. In particular, we have the following lemma.
\begin{lemma} \label{lem:couner-example}
For the family of instances \eqref{eq:lot-sizing-instance}, the optimal adjustable solution is $ z_{\sf AR}({\cal U}) = 0$ and the optimal affine solution is $ z_{\sf Aff}({\cal U}) = m/2-1$. In particular the gap between affine and adjustable policies is unbounded.
\end{lemma}
The proof of Lemma \ref{lem:couner-example} is presented in Appendix  \ref{apx-lem:couner-example}. Lemma \ref{lem:couner-example} shows that the assumption on the non-negativity of the recourse matrix $\mb B$ is necessary and crucial to obtain the theoretical bounds in  Table \ref{tab:results}. Relaxing this assumption can result in an unbounded gap. It is an interesting question to develop approximation algorithms and policies for two-stage robust problem with provable theoretical guarantees when the recourse matrix has negative components, or in particular is a network matrix.
}

\section{Conclusion.}

In this paper, we consider a two-stage adjustable robust problem~\eqref{eq:ar} under budget of uncertainty set and show that surprisingly affine policies give an optimal approximation for problem matching the hardness of approximation. We also present strong theoretical performance for more general intersection of budget of uncertainty sets that significantly improve over the worst case bound of $\Theta(\sqrt m)$. Budgeted uncertainty sets are an important class of uncertainty sets that are widely used in practice. Therefore, our improved bounds for affine policies for this class of uncertainty sets also address the stark contrast between the near-optimal empirical performance and the worst case performance of 
$\Theta(\sqrt m)$. 

Furthermore, our analysis shows that there exists a near-optimal affine solution satisfying a nice structural property where the scenarios are partitioned into {\em inexpensive} and {\em expensive} based on a threshold and the affine solution covers \red{only the inexpensive components (or a fraction of them) using a linear solution and remaining demand using a static solution.} This structure is closely related to threshold policies that are widely used in many applications, and allows us to design an alternate algorithm for computing near-optimal affine solutions for budget of uncertainty sets that is significantly faster than solving a large LP. This structural property might be of independent interest for other applications \red{and could provide insights to design more general policies that work well in settings where affine policies could be highly sub-optimal.}

\section{Acknowledgments.} The research is partially supported by NSF grants CMMI 1351838, CMMI 1636046 and DARPA Lagrange grant.

%


%
%
%

\bibliographystyle{informs2014} 
\bibliography{robust}

\begin{APPENDICES}

\section{Proof of Claims in Section \ref{section:claims}} \label{appendix:lem:structural-resulto}
\

\

\proof{Proof of Claim~\ref{claim:3.1}.}
\begin{align*}
z({\cal W})   &= \min_{\mb y \geq \mb 0} \left\{ \mb d^T \mb y \; \bigg\vert \; \mb B \mb y \geq \sum_{ i \in {\cal W}} \mb e_i  \right\} \\
&= \min_{\mb y \geq \mb 0 } \left\{ \mb d^T \mb y \; \bigg\vert \; \sum_{j=1}^n B_{ij} y_j \geq 1, \; \forall i \in {\cal W} \right\} \\
&= \min_{\mb y \geq \mb 0 } \left\{ \mb d^T \mb y \; \bigg\vert \; \sum_{j=1}^n  \frac{w_i B_{ij}}{ \max_{k\in {\cal W}} (w_k B_{kj})} \cdot   \max_{k \in {\cal W}} (w_k B_{kj}) \cdot y_j \geq w_i, \; \forall i \in {\cal W} \right\} \\
&= \min_{\mb y \geq \mb 0 } \left\{ \mb d^T \mb y \; \bigg\vert \; \sum_{j=1}^n \hat{B}_{ij} \cdot \max_{k \in {\cal W}}  (w_k B_{kj}) \cdot  y_j\geq w_i, \; \forall i \in {\cal W} \right\} \\
&= \min_{\mb y \geq \mb 0 } \left\{ \sum_{j=1}^n  \frac{d_j}{\max_{k \in {\cal W}} (w_k B_{kj})} \cdot  y_j  \; \bigg\vert \; \sum_{j=1}^n \hat{B}_{ij}   y_j\geq w_i, \; \forall i \in {\cal W} \right\} \\
&=  \min_{\mb y \geq \mb 0 } \left\{ \eta \gamma \sum_{j=1}^n  \hat{d}_j y_j  \; \bigg\vert \; \sum_{j=1}^n \hat{B}_{ij}   y_j\geq w_i, \; \forall i \in {\cal W} \right\} \\
&= \eta \gamma \cdot    \hat z({\cal W})  .
\end{align*}
\hfill
\Halmos
\endproof

\vspace{2mm}

\proof{Proof of Claim~\ref{claim:3.2}.}
Consider $ j \in [n]$, by the feasibility of the solution $ \mb{z^*}$ we have,
$$ \sum_{i \in {\cal J}} \hat{B}_{ij} \cdot z_i^* \leq\hat{d}_j ,$$
and consequently
$$\sum_{i \in {\cal J}} \hat{B}_{ij} \frac{2 \lfloor{z_i^*}\rfloor}{\eta}  \leq \sum_{i \in {\cal J}} \hat{B}_{ij} \frac{2 z_i^*}{\eta} \leq   \frac{ 2\hat{ d}_j}{\eta} . $$
Hence,
\begin{align*}
\mathbb{P} \left(  \sum_{i \in {\cal J}} \hat{B}_{ij}  \cdot \frac{ 2Z_i}{\eta}      > \hat{d}_j  \right) &=\mathbb{P} \left(  \sum_{i \in {\cal J}} \hat{B}_{ij} \cdot \frac{ 2 ( \lfloor{z_i^*}\rfloor+ {\xi}_i)}{\eta}     >  \hat{d}_j  \right)  \\
&\leq \mathbb{P} \left(  \sum_{i \in {\cal J}} \hat{B}_{ij} \frac{2 {\xi}_i}{\eta}     > (1- \frac{2}{\eta}) \hat{d}_j  \right)  \\
&=\mathbb{P} \left(  \sum_{i \in {\cal J}} \hat{B}_{ij}  {\xi}_i    > (\frac{\eta}{2}- 1) \hat{d}_j  \right).
\end{align*}
Now, we apply the Chernoff inequality in Lemma \ref{lem:chernof1} with $\delta =\frac{\eta}{2} -2$ and $\Xi=  \sum_{i \in {\cal J}} \hat{B}_{ij}  {\xi}_i$. Note that 
$ \delta = 2 \cdot \frac{\log n }{\log \log n}-2 >0$ for sufficiently large $n$.  Moreover, we have for all $ i \in {\cal J}, j \in {\cal J}$, $\hat{B}_{ij} \in [0,1]$ and
$$\mathbb{E} \left( \sum_{i \in {\cal J}} \hat{B}_{ij} {\xi}_i \right)  =  \sum_{i \in {\cal J}} \hat{B}_{ij} (z_i^*- \lfloor{z_i^*}\rfloor) \leq \sum_{i \in {\cal J}} \hat{B}_{ij} z_i^*  \leq \hat{d}_j.$$
Therefore the Chernoff bound gives,
$$
\mathbb{P} \left(  \sum_{i \in {\cal J}} \hat{B}_{ij}  {\xi}_i    > (\frac{\eta}{2}- 1) \hat{d}_j  \right)  \leq  \left( \frac{e^{\frac{\eta}{2} -1}}{ (\frac{\eta}{2}-2)^{\frac{\eta}{2} -2}} \right)^{\hat{d_j}}.$$
Recall $ \eta = 4 \frac{\log n}{\log \log n}$. Hence the RHS is equivalent to 
\begin{align*}
\left( \frac{e^{\frac{\eta}{2} -1}}{ (\frac{\eta}{2}-2)^{\frac{\eta}{2} -2}} \right)^{\hat{d_j}} = O \left( \frac{e^{\frac{\eta \hat{d_j}}{2}}}{ (\frac{\eta}{2})^{\frac{\hat{d_j} \eta}{2} }} \right) & = 
O  \left( \exp \left( \hat{d_j}\frac{\eta }{2} \left( 1 -  \log (\frac{\eta}{2}) \right) \right) \right) \\ 
 & = O  \left( \exp \left( -\hat{d_j}\frac{\eta}{2}  \log \eta \right) \right) \\
  & = O  \left( \exp \left( -2\hat{d_j} \frac{\log n}{\log \log n}  \log \left(4 \frac{\log n}{\log \log n} \right) \right) \right) \\
    & = O  \left( \exp \left( -2\hat{d_j}\log n \right) \right)  = O \left(\frac{1}{n^{2\hat{d_j}}} \right) \leq \frac{c}{n^2},
\end{align*}
for some constant $c$. The last inequality holds because $\hat{d}_j \geq 1$. Hence, 
$$\mathbb{P} \left(  \sum_{i \in {\cal J}} \hat{B}_{ij}  \cdot \frac{ 2Z_i}{\eta}      > \hat{d}_j  \right)  \leq \frac{c}{n^2}.$$
Therefore by a union bound we have,
$$ \mathbb{P} \left(  \sum_{i \in {\cal J}} \hat{B}_{ij} \frac{ 2 Z_i}{\eta}    > \hat{d}_j , \; \exists j \in [n] \right)  \leq \sum_{i=1}^n \mathbb{P} \left(  \sum_{i \in {\cal J}}\hat{B}_{ij} \frac{2 Z_i}{\eta}    > \hat{d}_j  \right)  \leq \frac{c}{n}.$$
Therefore,
$$ \mathbb{P} \left(  \sum_{i \in {\cal J}} \hat{B}_{ij} \frac{ 2 Z_i}{\eta}    \leq  \hat{d}_j , \; \forall j \in [n] \right)  \geq 1- \frac{c}{n}=1- O(\frac{1}{n}).$$
\hfill
\Halmos
\endproof

\vspace{2mm}

\proof{Proof of Claim~\ref{claim:3.3}.}
We have, 
$$ \sum_{ i\in {\cal J}} w_i  Z_i = \lambda +  \sum_{i \in {\cal J}}w_i {\xi}_i,$$
where
$$ \lambda =  \sum_{i \in {\cal J}} w_i \lfloor{z_i^*}\rfloor .$$
We apply the Chernoff inequality in Lemma \ref{lem:chernof2} with  $ \delta = \frac{1}{2}$ and $\Xi = \lambda +  \sum_{i \in {\cal J}}w_i {\xi}_i$. Note that $$\mathbb{E}(\Xi )=  \sum_{i \in {\cal J}}w_i  z_i^* > 1.$$ Hence,
\begin{align*}
\mathbb{P} \left(\sum_{i \in {\cal J}}w_i Z_i  > \frac{1}{2}\right) & =\mathbb{P} \left(\Xi >\frac{1}{2}\right)\\
&\geq \mathbb{P} \left(\Xi >\frac{\mathbb{E}(\Xi)}{2}\right)\\
& \geq 1 - e^{  -\frac{ \mathbb{E}(\Xi)}{8}}  \geq 1 - e^{  -\frac{ 1}{8}}.
\end{align*}
\hfill
\Halmos
\endproof

\section{Chernoff bounds}

  \begin{lemma}[\bf Chernoff Bound 1]\label{lem:chernof1}
 Let ${\xi}_1,{\xi}_2,\ldots,{\xi}_r$ be independent Bernoulli trials. Denote $\Xi = \sum_{i=1}^r \alpha_i {\xi}_i$ where $\alpha_1, \ldots, \alpha_r $ are reals in $[0,1]$. Let $s >0$ such that $\mathbb{E}(\Xi) \leq s$. Then for any $ \delta >0$,
 $$ \mathbb{P}( \Xi \geq (1+ \delta ) s ) \leq  \left(  \frac{e^{\delta}}{(1+\delta)^{1+\delta} }    \right)^s.$$
 \end{lemma}
The Chernoff bound in Lemma \ref{lem:chernof1} is a slight variant of the Raghavan-Spencer inequality (Theorem 1 in \cite{raghavan1986probabilistic}).  The proof is along the same lines as in  \cite{raghavan1986probabilistic}. For completeness, we are providing it below.

\proof{Proof.}
From Markov's inequality we have for all $t>0$,
 $$ \mathbb{P}( \Xi \geq (1+ \delta ) s ) = \mathbb{P}( e^{t\Xi} \geq e^{t(1+\delta)s} ) \leq \frac{\mathbb{E}(e^{t\Xi})}{e^{t(1+\delta)s}}.$$
 Denote $p_i$ the parameter of the Bernoulli ${\xi}_i$. By independence, we have
 $$  \mathbb{E}(e^{t\Xi})= \prod_{i=1}^r \mathbb{E} (e^{t \alpha_i {\xi}_i}) = \prod_{i=1}^r \left( p_i  e^{t \alpha_i }+1-p_i \right)
 \leq \prod_{i=1}^r \exp \left( p_i  (e^{t \alpha_i }-1) \right)$$
 where the inequality holds because $1+x \leq e^x$ for all $x \in \mathbb{R}$.
 By taking $t=\ln( 1+ \delta) >0$, the right hand side becomes
 \begin{align*}
 \prod_{i=1}^r \exp \left( p_i  ((1+\delta)^{ \alpha_i }-1) \right) &  \leq \prod_{i=1}^r \exp \left( p_i  \delta \alpha_i \right) \\
&  =\exp \left(   \delta  \cdot \mathbb{E}(\Xi ) \right) \leq e^{\delta s} ,\\
 \end{align*}
where the first inequality holds because $ (1+x)^{\alpha} \leq 1+ \alpha x$ for any $x \geq 0$ and $\alpha \in [0,1]$ and the second one because $s \geq \mathbb{E}(\Xi )= \sum_{i=1}^r \alpha_i p_i$. Hence, 
we have  $$  \mathbb{E}(e^{t\Xi}) \leq e^{\delta s}. $$
On the other hand, $$e^{t(1+\delta)s} = (1+\delta) ^{(1+\delta)s}.$$
Therefore,
 $$ \mathbb{P}( \Xi \geq (1+ \delta ) s ) \leq  \left(  \frac{e^{\delta s}}{(1+\delta) ^{(1+\delta)s}} \right) =\left(  \frac{e^{\delta}}{(1+\delta)^{1+\delta} }    \right)^s. $$
\hfill
\Halmos
\endproof

\begin{lemma} [\bf Chernoff Bound 2] \label{lem:chernof2}
 Let ${\xi}_1,{\xi}_2,\ldots,{\xi}_r$ be independent Bernoulli trials. Denote $$\Xi = \lambda+ \sum_{i=1}^r \alpha_i {\xi}_i$$ where $\lambda \in \mathbb{R}_+$ and $\alpha_1, \ldots, \alpha_r $ are reals in $(0,1]$. Denote $\mu=\mathbb{E}(\Xi )$. Then for any $ 0<\delta < 1,$  
$$ \mathbb{P} ( \Xi > (1- \delta ) \mu ) > 1 - e^{  \frac{-{\delta}^2 \mu}{2}}.$$
\end{lemma}
This is a slight variant of the  lower tail Chernoff bound \cite{chernoff1952measure}. The proof is along the same lines as in \cite{chernoff1952measure}.  For completeness, we are providing it below.

\proof{Proof.}
We show equivalently that 
$$ \mathbb{P} ( \Xi \leq (1- \delta ) \mu ) \leq  e^{  \frac{-{\delta}^2 \mu}{2}}.$$
From Markov's inequality we have for all $t< 0$,
 $$ \mathbb{P}( \Xi \leq (1- \delta ) \mu ) = \mathbb{P}( e^{t\Xi} \geq e^{t(1-\delta)\mu} ) \leq \frac{\mathbb{E}(e^{t\Xi})}{e^{t(1-\delta) \mu}}.$$
 Denote $p_i$ the parameter of the Bernoulli ${\xi}_i$. By independence, we have
 $$  \mathbb{E}(e^{t\Xi})=  e^{t \lambda}    \prod_{i=1}^r \mathbb{E} (e^{t \alpha_i {\xi}_i}) = e^{t \lambda} \prod_{i=1}^r \left( p_i  e^{t \alpha_i }+1-p_i \right)
 \leq e^{t \lambda}  \prod_{i=1}^r \exp \left( p_i  (e^{t \alpha_i }-1) \right) ,$$
 where the inequality holds because $1+x \leq e^x$ for all $x \in \mathbb{R}$.
We take $t=\ln( 1- \delta) < 0$. We have $ t \leq - \delta$, hence
$$ e^{t \lambda} \leq e^{- \delta \lambda} .$$
Moreover,
 \begin{align*}
\prod_{i=1}^r \exp \left( p_i  (e^{ t\alpha_i }-1) \right)   =\prod_{i=1}^r \exp \left( p_i  ((1-\delta)^{ \alpha_i }-1) \right) 
\leq \prod_{i=1}^r \exp \left( -p_i  \delta \alpha_i \right),
 \end{align*}
where the  inequality holds because $ (1-x)^{\alpha} \leq 1- \alpha x$ for any $ 0 < x <1$ and $\alpha \in [0,1]$. Therefore,
 $$  \mathbb{E}(e^{t\Xi}) \leq      e^{- \delta \lambda}       \prod_{i=1}^r \exp \left( -p_i  \delta \alpha_i \right)   =          e^{-\delta \mu}. $$
On the other hand, $$e^{t(1-\delta) \mu} = (1-\delta) ^{(1-\delta)\mu}.$$
Therefore,
 $$ \mathbb{P}( \Xi \leq (1- \delta ) \mu ) \leq  \left(  \frac{e^{-\delta \mu}}{(1-\delta) ^{(1-\delta)\mu}} \right) =\left(  \frac{e^{-\delta}}{(1-\delta)^{1-\delta} }    \right)^{\mu}. $$
Finally, we have for any $0<\delta < 1$,
$$ \ln (1- \delta ) \geq  - \delta + \frac{\delta^2}{2}$$
which implies 
$$ (1- \delta) \cdot \ln (1- \delta ) \geq  - \delta + \frac{\delta^2}{2}$$
and consequently
$$\left(  \frac{e^{-\delta}}{(1-\delta)^{1-\delta} }    \right)^{\mu} \leq e^{  \frac{-{\delta}^2 \mu}{2}}. $$
\hfill
\Halmos
\endproof

\section{Proof of Lemma \ref{lem:LP-form}} \label{apx-LP}
\

\

The affine problem \eqref{eq:aff} has the following epigraph formulation
 \begin{equation*}
\begin{aligned}
z_{\sf Aff}({\cal U}) = \min \; & \mb{c}^T \mb{x} + z\\
& z \geq   \mb{d}^T\left( \mb{P}\mb{h}+\mb{q}\right),  \; \; \; \forall \mb h \in {\cal U} \\
& \mb{A}\mb{x} + \mb{B}\left( \mb{P}\mb{h}+\mb{q}\right) \; \geq \; \mb{h},   \; \; \; \forall \mb h \in {\cal U} \\
& \mb{P}\mb{h}+\mb{q} \; \geq \; \mb{0},   \; \; \; \forall \mb h \in {\cal U} \\
& \red{\mb{x}  \in  {\cal X} }, \; \mb P \in \mathbb{R}^{n \times m}, \; \mb q \in \mathbb{R}^n, \; z \in \mathbb{R}.   \\
\end{aligned}
\end{equation*}
We use standard duality techniques to derive formulation \eqref{eq:aff:LP}. The first constraint is equivalent to 
$$ z - \mb d^T \mb q \geq \underset {\mb h \geq \mb 0 }{\underset{  \mb R  \mb h \leq \mb r  }\max } \; \mb{d}^T\mb{P}\mb{h} .$$
By taking the dual of the maximization problem, the constraint is equivalent to
$$ z - \mb d^T \mb q \geq \underset {\mb v \geq \mb 0 }{\underset{   \mb R^T \mb v \geq \mb P^T \mb d  }\min } \;      \mb{r}^T\mb v .$$
We can then drop the min and introduce $\mb v $ as a variable, hence we obtain the following linear constraints
$$ z- \mb d^T \mb q \geq \mb r^T \mb v $$
$$ \mb R^T \mb v \geq \mb P^T \mb d $$
$$ \mb{v}   \in  {\mathbb R}^{L}_+.$$
We use the same technique for the second sets of constraints, i.e.,
$$ \mb{A}\mb{x} + \mb{B}   \mb q   \ \; \geq \;    \underset {\mb h \geq \mb 0 }{\underset{  \mb R  \mb h \leq \mb r  }\max } \;         \mb{h} (  \mb I_m - \mb B \mb P )  .$$
By taking the dual of the maximization problem for each row and dropping the $\min$ we get the following compact formulation of these constraints
$$ \mb A \mb x + \mb B \mb q \geq \mb V^T \mb r $$
$$ \mb R^T \mb V \geq \mb I_m - \mb B \mb P $$
$$  \mb{V}   \in  {\mathbb R}^{L \times m}_+ .$$
Similarly, the last constraint 
$$ \mb{q} \; \geq \;\underset {\mb h \geq \mb 0 }{\underset{  \mb R  \mb h \leq \mb r  }\max } \;      -\mb{P}\mb{h} ,  $$
is equivalent to 
$$ \mb q \geq \mb U^T \mb r $$
$$ \mb U^T \mb R + \mb P \geq \mb 0 $$
$$\mb{U}   \in  {\mathbb R}^{L \times n}_+ . $$
Putting all together, we get the formulation  \eqref{eq:aff:LP}.

\section{Proof of Lemma \ref{lem:couner-example}} \label{apx-lem:couner-example}

\

\
First let us show that $z_{\sf AR}= 0$.
We consider the following solution for the  adjustable problem  
$$
x_i= \left\{
    \begin{array}{ll}
        1 & \mbox{if } i \in J_1 \\
        0 & \mbox{if }  i \in J_2.
    \end{array}
\right.
$$
Consider a scenario $\mb h $ that is an extreme point of $  {\cal U}$. In particular, we have  $h_i \in \{0,1\}^m$ for all $i$ and $\sum_{i=1}^m h_i = m/2$. Consider $\tilde{J}_1(\mb h )$ and $\tilde{J}_2(\mb h )$   respectively subsets of $ J_1 $  and $  J_2 $ in  which the demand is realized. In particular,
$$ \tilde{J}_1( \mb h)= \{ i \in J_1 \; \vert \; h_i =1 \} $$
$$ \tilde{J}_2 ( \mb h)= \{ i \in J_2 \; \vert \; h_i =1 \}.$$
We have, $$ \vert \tilde{J}_1 (\mb h) \vert +  \vert \tilde{J}_2 (\mb h ) \vert  = m/2 .$$  
Note that the first stage solution $\mb x$ covers the demand of the nodes in $ \tilde{J}_1 ( \mb h)$ because $x_i=1$ for all $i \in J_1$. On the other side, we cover demand  in $ \tilde{J}_2 ( \mb h)$ by sending inventory from $J_1 \setminus  \tilde{J}_1(\mb h) $ to $ \tilde{J}_2 ( \mb h) $ in the second stage. This is possible because
$$ \vert J_1 \setminus  \tilde{J}_1(\mb h) \vert = m/2 - \vert \tilde{J}_1(\mb h) \vert  =  \vert \tilde{J}_2 (\mb h) \vert . $$ The cost of sending inventory from $J_1 \setminus  \tilde{J}_1(\mb h)$  to $ \tilde{J}_2 (\mb h )$ is 0 because all directed distances from $J_1$ to $J_2$ are  zero.  In particular, we consider a matching $M$  from $J_1 \setminus  \tilde{J}_1(\mb h)$  to $ \tilde{J}_2 (\mb h )$. We define the following second stage solution

$$
y_{ij}(\mb h)= \left\{
    \begin{array}{ll}
        1 & \mbox{if } (i,j) \in M \\
        0 & \mbox{otherwise}.
    \end{array}
\right.
$$
We have $\mb x, \mb y( \mb h)$ is feasible for the adjustable problem and its corresponding cost is 0. Therefore, $ z_{\sf AR}({\cal U}) = 0$.

Now, let us show that $ z_{\sf Aff}({\cal U}) = m/2-1$. Given the distances in the graph, flow can be sent only from $J_1$ to $J_2$. In particular, we can rewrite the covering constraints of the problem as follows
\begin{equation} \label{morocco1}
\forall j \in J_2 \qquad x_j + \sum_{i \in J_1} y_{ij} (\mb h) \geq h_j
\end{equation}
\begin{equation} \label{morocco2}
\forall i \in J_1 \qquad x_i- \sum_{j \in J_2} y_{ij} (\mb h) \geq h_i. 
\end{equation}

Consider $\mb x$ and $\mb y( \mb h ) = \mb P \mb h + \mb q$ a feasible  affine solution. The number of rows of $\mb P$ is the number of edges in the graph. The number of columns of $\mb P$ is $m$ which is the total number of nodes. In particular,
$$     \forall i \in J_1,  \;   \forall j \in J_2 \qquad      y_{ij}( \mb h ) = \sum_{\ell =1}^m P_{(i,j),\ell} h_{\ell} + q_{ij},$$
where $P_{(i,j),\ell} $ denotes the component of $\mb P$ corresponding to edge $(i,j)$ and node $\ell $. Nodes of $J_1$ should be covered in the first stage because flow can not be sent to $J_1$ in the second stage. Moreover, $J_1$ is  covered at a zero cost.  In particular we have,
$$ x_i = 1 \qquad \forall i \in J_1.$$ 
Let fix $i \in J_1$. Consider $\mb h \in {\cal U}$ such that $h_i=1$. From \eqref{morocco2}, we have, 
$$ x_i \geq \sum_{j \in J_2} y_{ij}( \mb h ) + 1.$$
Moreover, we know that $ x_i \leq 1$ and $ y_{ij}( \mb h ) \geq 0$ for all $j \in J_2$. This implies that for all $j \in J_2$,
$$ y_{ij}( \mb h ) =0,$$
which is equivalent to 
$$ P_{(i,j),i}+ q_{ij}+ \sum_{\ell =1, \ell \neq i}^m P_{(i,j),\ell} h_{\ell} =0$$
for any $\mb h \in {\cal U}$ such that $h_i=1$. Therefore, for any $i \in J_1$, for any $j \in J_2$,
$$ P_{(i,j),i}+ q_{ij} =0,$$
and 
\begin{equation} \label{eq:imo}
 P_{(i,j),\ell} =0  \qquad \forall {\ell \neq i}  .
\end{equation}
Now fix $\ell \in J_1$ and $j \in J_2$. Consider the following scenario

$$
 \left\{
    \begin{array}{ll}
        h_{\ell} =0  \\
        h_i = 1 \qquad \forall i \in J_1 \setminus \{\ell\} \\
        h_j = 1 \\
        h_k= 0  \qquad \forall k \in J_2 \setminus \{j\}.
    \end{array}
\right.
$$
From \eqref{morocco1}, we have for all $j \in J_2$,
$$ x_j + \sum_{i \in J_1} y_{ij} ( \mb h) \geq 1.$$
Since there is a unit demand at node $i \in J_1$ for any $i \neq \ell$, we can not send inventory from $i$ to any node in $J_2$. In particular, $y_{ij}(\mb h)= 0 $ for any $i \neq \ell$. This implies from the last inequality that
$$ x_j + y_{\ell j} ( \mb h) \geq 1,$$
which is equivalent to 
$$ x_j + P_{(\ell,j),j} + \sum_{k \in J_1, k\neq \ell}P_{(\ell,j),k}  + q_{\ell j} \geq 1.$$
Putting it together with \eqref{eq:imo}, we get for any $j \in J_1$ and for any $\ell \in J_2$, we have
$$ x_j + q_{\ell j} \geq 1.$$
Therefore,
$$ \sum_{j \in J_2} x_j + \sum_{j \in J_2} q_{\ell j}  \geq m/2.$$
Now consider \eqref{morocco2} for $ \mb h =\mb 0$. We have for $ \ell \in J_1$,
$$ x_{\ell} - \sum_{j \in J_2}  q_{\ell j} \geq 0.$$
Moreover since $x_{\ell} \leq 1$ , we conclude that,
$$ \sum_{j \in J_2}  q_{\ell j}  \leq 1.$$
Therefore, 
$$ \sum_{j \in J_2} x_j   \geq m/2 -1.$$
Finally,
$$ z_{\sf Aff}  \geq   \sum_{j=1}^m  c_j x_j   =  \sum_{j \in J_2}  x_j  \geq m/2 -1.$$
Now we consider the following affine solution
$$
 \left\{
    \begin{array}{ll}
        x_i= 1 \qquad \forall i=1, \ldots, m-1  \\
        x_m= 0
    \end{array}
\right.
$$

$$
\forall {\mb h} \in {\cal U},  \forall i \in J_1, \left\{
    \begin{array}{ll}
        y_{ij}(\mb h)= 0  \qquad  \forall j \in J_2 \setminus \{m \}  \\
        y_{im}(\mb h)= -h_i +1.
    \end{array}
\right.
$$
The above affine solution is feasible for the adjustable problem. In fact, The capacity constraints are verified on $\mb x$. The non-negativity constraints are verified on $\mb y( \mb h) $. Demand is covered at each node by the first stage solution $\mb x$  except a node $m$. Demand at node $m$ is covered in the second stage because
 $$ \sum_{i \in J_1} y_{im}(\mb h) = \sum_{i \in J_1} (1 -h_i) = m/2 - \sum_{i \in J_1}  h_i  \geq h_m,$$ 
 where the inequality follows from the definition of $\cal U$.  The cost of this affine solution is $m/2 -1$. Therefore,
 $$ z_{\sf Aff} ( {\cal U})= m/2 -1.$$

\end{APPENDICES}





\end{document}